\documentclass[11pt,reqno]{amsart}
\usepackage{pdfpages}
\usepackage[english]{babel}
\usepackage[utf8]{inputenc} 
\usepackage[T1,OT1]{fontenc} 

\usepackage{amssymb} 
\usepackage{amsmath} 
\usepackage{amsthm} 
\usepackage{graphicx,subfigure}	
\usepackage{fullpage} 
\usepackage{epsfig,color}
\usepackage{enumerate,verbatim}
\usepackage[toc,page]{appendix}
\usepackage[square,numbers]{natbib}
\usepackage[final]{showlabels}
\usepackage{enumitem}

\usepackage{float} 

\usepackage{tikz}
\usetikzlibrary{shapes}
\usepackage{pgfplots}
\usepackage{todonotes}

\definecolor{green_mermoz}{rgb}{0,0.5,0}
\usepackage{hyperref}
\usepackage{fontawesome}

\definecolor{b}{rgb}{0.00000,0.44700,0.74100}%
\definecolor{o}{rgb}{0.85000,0.32500,0.09800}%
\definecolor{y}{rgb}{0.92900,0.69400,0.12500}%
\definecolor{g}{rgb}{0.46600,0.67400,0.18800}%
\definecolor{p}{rgb}{0.49400,0.18400,0.55600}%
\definecolor{inrae}{RGB}{0,163,166}
\definecolor{k}{rgb}{0,0,0}%
\definecolor{gray}{rgb}{0.75,0.75,0.75}%
\definecolor{pink}{rgb}{1,0.7529,0.7961}%
\definecolor{r}{rgb}{0.6350,0.0780,0.1840}%
\definecolor{cyan}{rgb}{0.3010,0.7450,0.9330}%
\hypersetup{colorlinks = true, linkcolor=b, citecolor=o}
\bibpunct{\textcolor{o}{[}}{\textcolor{o}{]}}{,}{n}{}{;}

\usepackage[noabbrev]{cleveref}
\creflabelformat{equation}{#2(#1)#3}
\crefformat{equation}{#2(#1)#3}
\crefformat{equations}{#2(#1)#3}
\crefname{thm}{Theorem}{Theorems}
\crefname{ass}{Assumption}{Assumptions}
\crefname{problem}{Problem}{Problems}
\crefname{prop}{Proposition}{Propositions}
\crefname{lem}{Lemma}{Lemmas}
\creflabelformat{problem}{#2(#1)#3}

\definecolor{k4}{rgb}{0.8,0.8,0.8}
\definecolor{k3}{rgb}{0.6,0.6,0.6}
\definecolor{k2}{rgb}{0.4,0.4,0.4}
\definecolor{k1}{rgb}{0.2,0.2,0.2}

\newcommand{\N}{\mathbb{N}}

\newcommand{\R}{\mathbb{R}}

\newcommand{\p}{\partial}

\newcommand{\abs}[1]{\left\lvert#1\right\rvert}
\newcommand{\norm}[1]{\left\lVert#1\right\rVert}

\renewcommand{\epsilon}{\varepsilon}
\renewcommand{\phi}{\varphi}

\DeclareMathOperator*{\argmin}{argmin}   

\DeclareMathOperator{\zs}{z_\ast}

\renewcommand{\leq}{\leqslant}
\renewcommand{\geq}{\geqslant}

\renewcommand{\ge}{\geqslant}

\renewcommand{\tilde}{\widetilde}

\newcommand{\Var}{\mathrm{Var}}

\newcommand{\e}{\varepsilon}

\newcommand{\dint}{\displaystyle \int}

\newtheorem{thm}{Theorem} 
\newtheorem{prop}[thm]{Proposition} 
\newtheorem{cor}[thm]{Corollary}

\numberwithin{equation}{section}
\numberwithin{thm}{section}

\def\ds{\displaystyle}
\usepackage{epstopdf}

\title{Ancestral lineages in mutation selection equilibria with moving optimum}

\author{Raphaël Forien} 
\address{BioSP, INRAE, 84914, Avignon France} 
\email{raphel.forien@inrae.fr}
\author{Jimmy Garnier} 
\address{LAMA, UMR 5127 CNRS \& Univ. Savoie Mont-Blanc, Chambéry, France}
\email{jimmy.garnier@univ-smb.fr}
\author{Florian Patout} 
\address{BioSP, INRAE, 84914, Avignon France} 
\email{florian.patout@inrae.fr}
\date{\today}

\begin{document}
	
	\begin{abstract}
Many populations can somehow adapt  to rapid environmental changes. To understand this fast evolution, we investigate the genealogy of individuals inside those populations. More precisely, we use a deterministic model to describe the phenotypic density of a population under selection when the fitness optimum moves at constant speed. We study the inside dynamics of this population using the neutral fractions approach. We then define a Markov process characterizing the distribution of ancestral phenotypic lineages inside the equilibrium. This construction yields qualitative as well as quantitative properties on the phenotype of typical ancestors. In particular, we show that in asexual populations typical ancestors of present individuals carried traits much closer to the fitness optimum than most individuals alive at the same time. We also investigate more deeply the asymptotic regime of small mutation effects. In this regime, we obtain an explicit formula for the typical ancestral lineage using the description of the solutions of Hamilton Jacobi equation as a minimizer of an optimization problem. In addition, we compare our deterministic results on lineages with the lineages of stochastic models.

	\end{abstract}
	\maketitle	
	

	\section{Introduction}
	\subsection{Model description} \label{subsec:model}
	Numerous studies have reported rapid evolution in populations facing brutal environmental changes, such as climate change, habitat alteration or drug treatment~\cite{BraHol06,Hai05,hoffmann2011climate,parmesan2006}.
	These populations are under sustained pressure to adapt to changing environments. For instance, pathogens like influenza continuously adapt to evade their host's immune system, causing chronic infections despite heavy immune responses~\cite{BedCob11}. 
	Such rapid evolution is associated to specific genealogies, and leaves distinctive footprints in the genetic structure of populations.
	In particular, in asexual populations continuously adapting to changing environment, lineages trace back to a small pool of highly fit ancestors~\cite{neher2013genealogies}.
	Although many theoretical studies have focused on the statistical properties of genealogies using  coalescent models~\cite{Berestycki_recentprogress,BruDer07,Kingman82}, little is known about the lineages dynamics~\cite{companion}.
	
	
	Over the past few decades, important theoretical progress has been made in predicting phenotypic evolution in a changing environment (reviewed in~\citep{kopp2014rapid}). Since the pioneer work of~\cite{BurLyn95,LynGabWoo91,LynLan93}, most theoretical approaches to adaptation in a changing environment rely on a quantitative phenotypic trait subject to stabilizing selection around some optimal phenotype, whose value is shifted continuously through time.
	
	In this work, we investigate the genealogy of a population adapting to a changing environment. Our  aim is to discover which traits contribute to the adaptation of future generations. We first focus on the dynamics of lineages using a deterministic model. 
	Adapting the inside dynamics methods of~\cite{garnier2012inside,roques2012allee}, we track the traits of the progeny of different neutral fractions of the total population adapting to the changing environment. 
	Then, using duality between partial differential equations and stochastic processes,  
	we look backwards in time to investigate the ancestors of the population. 
	To this end, we define a stochastic  ancestral process that  describes the traits of the past ancestors of an individual sampled uniformly from those with a given trait in the present population. 
{Our analysis goes back and forth between a forwards deterministic description of the population adapting to a changing environment with the dynamics of neutral fractions inside, and a backwards description of ancestral lineages in the form of stochastic processes.}
	
	We consider a population characterized by a phenotypic trait $x\in\R$. Its trait density $f(t,x)$ changes as a result of mutations and natural selection.
	Individuals in this population give birth at a constant rate $\beta$. Their offspring inherit the parents' trait with possibly a variation due to mutations described by the mutation kernel $K_\sigma$. This yields the reproduction operator $\mathcal{B}$:
	\begin{equation}\label{eq:asexual B}
	\mathcal{B}(f)(x) := \big(K_\sigma * f\big)(x) =  \dfrac1{\sigma} \int_\R  K\left(\dfrac{x-x'}\sigma\right)f(x')\, dx'
	\end{equation}
	This convolution term accounts for an asexual reproduction of individuals. 
	We assume that $K$ satisfies:
	\begin{align}\label{ass K}
	\dint_\R K(y) dy = \dint_\R y^2 K(y)dy =1, \quad\quad \exists \eta>0 \text{ s.t. } \dint_\R K(y) e^{\eta \abs{y}} dy < + \infty,
	\end{align}
	so that the  parameter  $\sigma^2$ corresponds to the variance of the mutation kernel. 
	We assume that  selection  acts through the intrinsic mortality rate $d(t,\cdot)$. 
	Environmental change is modeled by assuming that this function is of the form
	\begin{align*}
	d(t,x) = \mu(x-ct),
	\end{align*}
	for some $c>0$ which measures the speed of environmental change.
	We assume that $\mu$ is a  \emph{convex}  function such that
	\begin{align}\label{hyp growth m}
	\lim_{\abs{x} \to  + \infty} \mu(x) = + \infty.
	\end{align}
	Without loss of generality we suppose that $\mu$ admits a (global) minimum at $x=0$ (so that the \emph{optimal trait} is always $ct$). 
	Finally, we assume that the death rate due to competition is given by $\gamma \int_\R f(t,z)dz$, 
	a non-local logistic competition term where  $\gamma > 0$ is  the competition factor.
	
	Overall, the evolution of $f$ 
	is described by the following integro-differential equation (IDE):
	\begin{equation} \label{eq:f z} 
	\partial_{t} f(t,x) = \beta \mathcal{B}( f(t,\cdot) )(x)\ -  \left(\mu(x - ct) + \gamma \ds \int_\R f(t,x')dx' \right)  f(t,x),  \text{ for } \, t>0 ,  \, x\in\R.
	\end{equation}
	In order to describe a population keeping pace with the changing environment, we look for traveling pulse solutions of \cref{eq:f z} which are special solutions of the form
	\begin{align}
	f(t,x) = F(x-ct),
	\end{align}
	where the constant profile $F$ describes a ``mutation-selection equilibrium'',  in the frame moving  at speed $c$.  It satisfies
	\begin{equation}\label{eq:TW_c}
	\gamma \int_\R F(z')dz' F(z)- c \partial_{z} F(z)   = \beta \mathcal{B}( F )(z)\ - \mu(z) F(z), \ z\in\R,
	\end{equation}
	where  $z$ represents the variable in the \emph{moving frame} at speed $c$,
	\begin{align*}
	z := x-ct.
	\end{align*}
	The phenotypic trait $z$ in the moving frame corresponds to the difference  between the trait $x$ and the moving optimal trait $ct$.  The existence of a traveling pulse solution of \cref{eq:TW_c} can be derived from the existence, stated in \cite{hamelcoville2019}, of a spectral pair $(\lambda,F)$ solving
	\begin{equation}\label{eq:equilibrium c}
	\lambda F(z) - c \partial_{z} F(z)   = \beta \mathcal{B}( F )(z) - \mu(z) F(z), \ z\in\R.
	\end{equation}
	In the following, we always assume that $\lambda > 0$.
	In this case, since \cref{eq:equilibrium c} is a linear equation, we can choose $F > 0$ such that
	\begin{equation}\label{eq:rhob}
	\lambda=\gamma \int_\R F(z)dz,
	\end{equation}
	yielding a solution to \cref{eq:TW_c}.
	The profile $F$ is bell shaped and centered around a value $\zs$. 
	This dominant trait $z^*$ is negative due to the advection term in \cref{eq:TW_c} and represents the lag of adaptation in a changing environment (see \Cref{fig initial frac}).
	
	
	\subsection{Mathematical results}
	In order to describe the ancestral lineages in the population, we  first track the offspring of individuals. To do so, we investigate the inside dynamics of the traveling pulse $F$ using the approach of \emph{neutral fractions} (see \Cref{fig neutral frac} for a schematic representation).
	The idea is that  individuals are labeled, and   transmit their label to their offspring. Since individuals only differ by their label and their trait, each label $k \in \N$ corresponds to a neutral fraction of density $\upsilon^k$ inside the population with density $F$. Initially, we assume 
	\begin{equation}\label{init frac}
	F(z) = \sum_{k\geq 1}\upsilon^k(0,z), \hbox{ for all } \, z\in\R.
	\end{equation}
	Their dynamics in the moving frame is described by
	\begin{equation} \label{eq:f_sub} 
	\left\{ \begin{array}{l}
	\ds \partial_{t} \upsilon^k(t,z) - c \p_z \upsilon^k(t,z) +  \Big(\mu(z) + \gamma \int_\R{F(z')\,dz'} \Big)  \upsilon^k(t,z) = \beta \mathcal{B}( \upsilon^k(t,\cdot) )(z)\, , \text{ for } \, t>0 , \, z \in\R\\
	\upsilon^k(0,z) = \upsilon_0^k(z), \ z \in\R.
	\end{array}\right.
	\end{equation}
	In particular, mutation and selection act on fractions the same way  they do on the entire population $F$ (solution of \cref{eq:TW_c}).  
	Note that by linearity, the sum $\sum_{k \ge 1}\upsilon^k(t,z+ct)$ verifies  \cref{eq:f z} with $f(0,x) = \sum_{k \geq 1} \upsilon^k(0,x)$ (fixed frame).
	Therefore,  with  the corresponding initial assumptions \cref{init frac}, we see that the neutral fractions are, at all times, a subdivision of the equilibrium profile :
	\begin{align}\label{ass subpop}
	\forall t\geq 0, \, \forall x \in \R : \quad  F(x-ct) = \sum_{k \ge 1}\upsilon^k(t, x-ct).
	\end{align}
	This approach has been introduced in the context of reaction diffusion equations to understand the evolution of diversity inside traveling wave solutions (in their context,  the front solution plays a similar role to our equilibrium profile $F$ : \cite{garnier2012inside,hallatschek2008gene,hallatschek2010life,roques2012allee}. 
	
	Here, the dynamics of neutral fractions describes the evolution of the progeny of different subgroups of individuals in the population, depending on the initial distribution of the subgroups. 
	As we shall see below, the distribution of the neutral fractions at any given point $x$ of the trait space can be deduced from the initial distribution of neutral fractions and the distribution of the ancestors of the individuals carrying the trait $x$.
	Combining this forward approach with the characterization of stochastic processes, we are able to define a backward \emph{ancestral process} in the moving frame.
	
	\begin{thm}[Ancestral process]\label{theo edp lineages}
		Let $\mathcal{A}$ denote the linear operator defined by
		\begin{align}\label{def A}
		\mathcal{A} \psi := \frac{\beta}{F}  \Big[  K_\sigma \ast (F \psi) - (K_\sigma \ast F)\, \psi  \Big] + c \p_y \psi,
		\end{align}
		for all $\psi \in \mathcal{C}^1_b(\R_+ \times \R)$.
		\begin{enumerate}
			\item The operator $\mathcal{A}$ generates a Feller semigroup, $(\mathcal{M}_s, s\geq 0)$, defined on the set of real continuous and bounded functions, $\mathcal C_b(\R)$.
			We call the Markov process associated to the semigroup $(\mathcal{M}_s,s \geq 0)$, denoted by $(Y_s, s \geq 0)$, the \emph{ancestral} process.
			\item The semigroup $(\mathcal{M}_s, s \geq 0)$ is such that, for all $(t,z) \in \R_+ \times \R$, and all $k\in \N$:
			\begin{align}\label{def M}
			\upsilon^k(t,z) = F(z) \mathcal{M}_{t-s}\left(\frac{\upsilon^k(s,\cdot)}{F}\right)(z), \  \, \forall\, 0 \leq s \leq t.
			\end{align}
			In other words, it is the moment dual of $\upsilon^k/F$, in the sense of stochastic processes, see \cref{dual process}.
		\end{enumerate}	
	\end{thm}
	From the definition of $(Y_s, s \geq 0)$, we obtain the following equality linking the forward model of neutral fractions with the backward ancestral process: for all $(t,z) \in \R_+ \times \R$ and $k\in \N$:
	\begin{align}\label{dual process}
	\frac{\upsilon^k(t,z)}{F(z)} = \mathbb{E}_z \left[  \frac{\upsilon^k(s,Y_{t-s})}{F(Y_{t-s})}\right], \ 0 \leq s \leq t.
	\end{align}
	This equality means that the probability of sampling an individual of type $k$ from those with trait $ z $ at time $ t $ (left hand side) is equal to the probability of drawing an individual of type $k$  among individuals of trait $y$ in the past, at the time $t-s$ (right hand side). 
	The right hand side measures the relative contributions of individuals living at time $ t-s $ to those carrying trait $ z $ at time $t$.
	Since these contributions are given by the distribution of $Y_s$, and \cref{dual process} is satisfied for any choice of $(\upsilon^k(0,\cdot), k \geq 1)$, we can say that $Y_s$ describes the trait of the ancestor (alive at time $t-s$) of an individual sampled at random from those carrying trait $z$ at time $t$.
	For this reason, we call $Y_s$ the \emph{ancestral process}.
	
	Note that the generator of $(Y_s, s \geq 0)$ comprises two parts: a jump part and a drift part.
	The process $(Y_s, s \geq 0)$ jumps at rate $\frac{\beta}{F(Y_s)} K_\sigma \ast F(Y_s)$ to a location $z$ whose density is given by
	\begin{align*}
	\frac{K_\sigma(Y_s-z) F(z)}{K_\sigma \ast F(Y_s)}.
	\end{align*}
	Moreover, between jumps, $Y_s$ drifts to the right at constant speed $c$.


	Note that the initial time of the ancestral process $Y_s$ (``$s=0$''), corresponds to \textbf{any} time $t \in \R$. The trajectory of ancestors is independent of the time $t$ at which we sample individuals in the population $F$ because the population is at equilibrium in the moving frame. 
	
	Next, we  establish some properties of \emph{lineages} which corresponds to trajectories of the ancestral process. We first start with the long-time behavior of the ancestral process  as $s \to + \infty$ which informs us on the trait distribution of the ancestors far back in the past.
	
	\begin{prop}(Common ancestors distribution)\label{long time asymp}$ $\\
		When $s \to \infty$, the ancestral process $Y_s$ converges in law towards a random variable $Y_\infty$, which admits the following density:
		$$\frac{ F \phi}{\int F(y') \phi(y')dy'},$$ 
		where $\phi$ is the non-negative stationary solution of the dual problem:
		\begin{align}\label{def dual L}
		\left( \gamma \int_\R F \right) \phi + c\partial_z \phi=  \beta\Big( K_\sigma* \phi - \phi\Big)  +  \Big( \beta - \mu   \Big) \phi \end{align}
	\end{prop}
	
	\cref{long time asymp} hinges on an intermediary result stated in \cite{gabrielcloez2019}. The description of the density of  $Y_\infty$ relies on the long time behavior analysis of the neutral fractions. 
	It also provides interesting properties on the traveling pulse $F$. 
	Indeed, we prove that any fraction inside the traveling pulse converges to a positive proportion of the traveling pulse as time goes to $\infty$. Thus, the traveling pulse is \emph{pushed} in the sense that any fraction inside this pulse push it forward~\cite{garnier2012inside,roques2012allee}.
	
	However, the common ancestor distribution is not uniform, 
	some traits  are more represented among ancestors than in the overall population.  In particular, we show that ancestors are more likely to have a trait close to the optimal trait $z=0.$  More precisely, we show that $Y_\infty$ satisfies these additional properties:
	\begin{cor}\label{cor explicit asymp}
		Under the additional assumption that selection is symmetric :  $\mu(x) = \mu(\abs x)$, the distribution of $Y_\infty$ admits the following density: 
		$$z \mapsto \frac{ F(z) F(-z) }{\int F(y') F(-y')dy'},$$ 
		In particular, 
		\begin{align*}
		\mathbb E \left[\psi(-Y_\infty)\right] = \mathbb{E}\left[ \psi(Y_\infty) \right], \ \hbox{for any } \  \psi \in \mathcal{C}_b(\R) \text{ and } 	\mathbb E \left[ Y_\infty \right]  =  0. 
		\end{align*}
	\end{cor}
	\Cref{cor explicit asymp} states that on average the 
	ancestors of the present population share the optimal trait with respect to selection. However, the population density at the optimal trait $F(0)$ is low because the dominant trait of the population lags behind the optimal trait at a distance  $z^*$ (see~\Cref{fig initial frac}). Thus, we learn  an interesting feature of the adaptation phenomenon of our model: the ancestors of typical individuals in the present population were far from typical within the population at the time.
	The existence of (potentially very few) optimally fit individuals is thus  crucial to the survival of the whole population since they will be the (most likely) ancestors of the next generations. 

	\bigskip 
	To further understand the trajectories of the ancestral process $Y_s$, we focus on the asymptotic regime of small mutations, that is $\sigma \to 0$ with a time of order $s/\sigma$. Under this asymptotic regime, we are able to fully characterize the ancestral process $Y_s$. 
	Since time is scaled by $1/\sigma$, we also need to rescale the speed of change $c$:
	\begin{equation}\label{scaled c}
	c = \sigma c'.
	\end{equation}
	We thus let $F_\sigma$ denote the solution to \cref{eq:TW_c} with $c$ replaced by $\sigma c'$, and $(\mathcal{M}_{\sigma, s}, s \geq 0)$ the associated lineage semigroup.
	\begin{prop}[Small mutations regime]\label{regime adaptative dynamics}$ $\\
		Let $\mathcal M ^\sigma_s:= \mathcal M_{\sigma,s/\sigma}$. Then, for all times $s \in \R_+$, $z \in \R$ and $\psi \in \mathcal{C}_b(\R)$, 
		\begin{align}\label{conv semigroup}
		\mathcal M^\sigma _s \psi(z) \xrightarrow[\sigma \, \to \, 0]{} \psi(\Gamma_z(s)), \text{ locally uniformly in time, }
		\end{align}
		and, for any $ \varepsilon > 0 $ and $ T > 0 $,
		\begin{align} \label{cvg_probability}
		\lim_{\sigma \, \to \, 0} \mathbb{P}\bigg( \sup_{s \in [0,T]} |Y_s - \Gamma_z(s)| > \varepsilon \, \Big| \, Y_0 = z \bigg) = 0.
		\end{align}
		The limit process $\Gamma_z$ corresponds to the solution of the following ODE:
		\begin{equation}\label{edo Gamma}
		\left\{ \begin{array}{ll}
		\dot \Gamma_z(s) & = c' - \beta \, \p_p H \Big[ U' (\Gamma_z (s)) \Big], \quad s>0 \\
		\Gamma_z(0) & = z.
		\end{array}\right.
		\end{equation}
		where the function $U$ and the Hamiltonian $H$ are defined by 
		\begin{align}\label{limit U}
		U = \lim_{\sigma \to 0} (-\sigma \log F_\sigma) \ \hbox{ and } \ 	H(p):=  \int_\R  K(y) \exp (y p ) dy -1.
		\end{align}	
	\end{prop}
	The limit of $\sigma \log(F_\sigma)$ in \cref{limit U} holds locally uniformly, and the existence of $U$ results from~\cite{lorz2011dirac}. In particular the authors show that $U$ is a $\mathcal C^1(\R)$ function, solution of an Hamilton Jacobi equation (see \Cref{sec HJ derivation} for a detailed derivation of this equation). Under this particular regime, the ancestral process becomes a deterministic process defined by $\Gamma_z$.  
	This situation was expected because the cumulative variance in the ancestral process due to mutations is proportional to $\sigma^2 (t/\sigma) = \sigma t $ in the rescaled time scale, and thus vanishes for any fixed $t$.
	
	This asymptotic regime has been widely studied in evolutionary contexts since the pioneer work of~\cite{diekmann2005}. This regime provides a good approximation when either the effects of mutations are small or when selection is strong, and describes how in this regime a population concentrates around one or several traits. This approach is connected to large deviations theory \cite{champagnathenry} and  it generally involves Hamilton-Jacobi equations~\cite{barles2009,lorz2011dirac}. Here, we show that the asymptotic equations keeps some trace of the history of ancestors, and thus  provide insights on the trajectories of the ancestral process that are the typical lineages.
	To the best of our knowledge, this aspect is new. In \Cref{sec HJ}, we detail how the Hamilton-Jacobi equation provides heuristics about typical lineages.  
	
	\subsection{Related works}
	Our results, stated for the non-local birth operator $\mathcal{B}$ defined by \cref{eq:asexual B} can be extended to the case where mutations are modeled by a Laplace operator ($\mathcal B = \Delta$). Some specific additional results are provided in   \Cref{sec diffusive approx}.
	In particular the trajectories of ancestral lineages follow an explicit SDE, see \cref{prop diff}. 
	
	The same diffusion process was obtained in \cite{companion} as the limit, as the number of individuals tends to infinity, of the trajectory of the traits along an individual's ancestral lineage in an individual-based stochastic population model.
	More precisely, the authors considered a stochastic model describing the evolution of a population of individuals following the rules described in \Cref{subsec:model}, with a finite number of individuals $N$ tending to infinity.
	They showed, using the historical process associated to their population model and a branching process approximation along with techniques initiated by \cite{marguetbranching}, that the lineage of an individual sampled uniformly from those alive at some time $T$ converges in distribution to the trajectory of a solution to the SDE \cref{SDE_Y_diffusive}, while at the same time the renormalized population process converges to a measure with density given by $F$, solution of \cref{eq:TW_c_diff}.
	Their result can be seen as a microscopic justification of our genealogical interpretation of the duality relation \cref{dual process}, in the case of a diffusive mutation operator. Recently, \cite{bansaye2021spine} used closely related probabilistic tools to those developed in \cite{companion} to describe ancestral lineages of  populations in a distinct model. 
	
	Mathematically, a large number of models describe the adaptation of populations in a \emph{steady} environment, (see e.g. \cite{diekmann2005,GilHamMarRoq19}). In the case of a linearly varying environment, there exist reaction diffusion models where a favorable region moves at a certain speed, \cite{BerDie09,BerFan18}, but it does not describe  an adaptation phenomenon contrary to \cite{AlfBer17}. Periodically fluctuating environment have also been studied, see \cite{FigMir18,FigMir19,LorChi15}. Recently, \cite{roques2020adaptation} proposed a methodology to  deal with general changing environments (linear, oscillating or stochastic) in the case of quadratic selection and a diffusive mutation operator.

	Our work has a close kinship with \cite{main}, where (mostly) formal analytical features measuring the dynamics of adaptation  are obtained, for  integro-differential models close to \cref{eq:f z}. The methodology is based upon asymptotic expansions in the  same regime as \Cref{regime adaptative dynamics}. It  also encapsulates 
	the case where the operator $\mathcal{B}$ describes sexual reproduction via the infinitesimal operator, see \cite{spectralsex,timespectralsex} for rigorous asymptotic treatment   (without environmental change) of this operator.

	\subsection{Ancestral lineages : an IDE point of view}\label{sec heur pde}
	The  fundamental solution associated to $\mathcal A$, the generator of the ancestral process $Y_s$, solves the following linear IDE: 
	\begin{equation}\label{pde w}
	\left\{	\begin{array}{ll} \p_s w^z(s,y) & = \beta \ds \int_{\R} K (h) \frac{F(y + \sigma h)}{F(y)}\Big( w^z(s,y+\sigma h)-w^z
	(s,y) \Big)dh + c \p_y w^z(s,y), \\
	w^z(0,y) & = \delta(z-y).
	\end{array}\right.
	\end{equation}
	This can be interpreted as a IDE describing the dynamics  of the trait $y$  of the ancestors whose descendants  have reached the trait $z$ in the population.

	Formally, the discrete  fractions label $k$ in \cref{eq:f_sub} is replaced by a continuum of neutral alleles. 
	Each neutral fraction then corresponds to the progeny of a single ancestral trait in the population.
	Indeed, let us define for any $y\in \R$,   $\upsilon^y$  as follows: 
	\begin{align*}
	\upsilon^y(t-s,z) = F(y) w^z(s,y) \, \text{ for  any } 0 \leq s \leq t, \text{ and } z \in \R,
	\end{align*}
	where $w^z$ is the fundamental solution  in \cref{pde w}.
	Let $\mathcal{L}$ be the linear operator defined by
	\begin{equation}\label{def L}
	\mathcal{L}(\psi) = \beta\Big( K_\sigma \ast \psi - \psi\Big) + c\partial_z \psi - \Big(\mu -\overline{\mu}  \Big)  \psi,
	\end{equation}
	where (integrating \cref{eq:TW_c} to obtain the last equality)
	\begin{align*}
	\overline{\mu} = \beta - \lambda = \int_\R \mu(z) \dfrac{F(z)}{\int_\R F(z')dz'} dz.
	\end{align*}
	A simple computation, detailed in \Cref{sec proof theo edp}, guarantees that  $\upsilon^y$ solves:
	\begin{equation}\label{eq:f_sub_dirac}
	\left\{ \begin{array}{ll}
	\p_t  \upsilon^y(t,z) & = \mathcal L(\upsilon^y(t,\cdot))(z),\\
	\upsilon^y(0,z) & = \delta(z-y)F(z).
	\end{array}\right.
	\end{equation}
This corresponds to the equation~\cref{eq:f_sub} satisfied by the neutral fractions :
	\begin{align}\label{def L frac}
	\forall k\in \N, \quad 	\p_s \upsilon^k = \mathcal{L}(\upsilon^k).
	\end{align}
	As a consequence,   $\upsilon^y$ can be seen as a neutral fraction, with an neutral label that singles out one trait in the population at $t=0$. 
	Therefore, \textbf{the (stochastic, backwards) ancestral process defined in \cref{theo edp lineages} corresponds to  a continuum of (deterministic, forward) neutral fractions.} Numerical simulations highlight this in \Cref{fig lineage asex distrib}, with more details provided in \Cref{app:IBM_vs_PDE}.
	
	We finally propose an interpretation in terms of partial differential equations of \cref{regime adaptative dynamics}.  The fundamental solution $w_\sigma^z$, introduced in \cref{pde w},   converges in the sense of distributions, when $\sigma \to 0$, up to the acceleration of time,  towards $w_0^z$ the solution of:
	\begin{equation}\label{pde transport}
	\left\{	\begin{array}{ll}
	\p_s w_0^z(s,y)& = - \Big(c'- \beta \, \p_p H(\p_y U(y))  \Big) \p_y w_0^z(s,y),\\
	w_0^z(0,y) & = \delta(z-y).
	\end{array}\right.
	\end{equation}
	The integral flow of the transport equation \cref{pde transport} coincides with the ODE \cref{edo Gamma} solved by $\Gamma$. 
	In \Cref{sec HJ},  
	we recover 
	this  formula (established in \cref{gamma flow}), directly from an Hamilton-Jacobi equation (rigorously derived in  \cite{lorz2011dirac}).  It has somehow  its origin in the field of the Weak-KAM theory, independently of neutral fractions. Therefore, thanks to $\Gamma$, our work links  the Hamilton-Jacobi framework and the ancestral process, in the regime $\sigma \to 0$. 
	
	\bigskip
	The rest of this article is organized as follows. In \Cref{sec proof theo edp}, we prove \cref{theo edp lineages}, using classical semigroup theory. 
	Then,  \Cref{sec:long_time_lineage} is devoted to the proof of \cref{long time asymp}. 
	In \Cref{sec diffusive approx}, we consider the  case of a difffusive mutation operator,  for which explicit computations are possible.
	Finally, we prove the \cref{regime adaptative dynamics} in \Cref{sec adap}. We also provide a discussion on how to obtain the ODE from heuristics on the Hamilton-Jacobi equation and a comparison of theoretical results with stochastic and deterministic simulations. 
	In the appendices, we provide additional discussions on the link between lineages and the methodology introduced in \cite{main} under the regime of small mutations. Finally, we detail the numerical methods to keep track of lineages with individual based simulations and the comparison with the deterministic model.
	
	
	\section*{Acknowledgements}
	The authors acknowledge Vincent Calvez for introducing them to the problem,  helping to formulate it and for his support during this work.
	We also warmly thank Jérôme Coville and Lionel Roques for their helpful comments. 
	This project has received funding from the European Research Council (ERC) under the European Union’s Horizon 2020 research and innovation programme (grant agreement No 639638) and from the French Agence Nationale de la Recherche (ANR-18-CE45-0019 "RESISTE" and ANR-16-CE02-0009 "GLOBNETS").
	R. F. was supported in part by the Chaire Modélisation Mathématiques et Biodiversité (École Polytechnique, Muséum national d'Histoire naturelle, Fondation de l'École Polytechnique, VEOLIA Environnement).

	\section{Link between lineages and fractions: proof of \cref{theo edp lineages}}\label{sec proof theo edp}
	(1) The existence of the  Feller semi-group holds true because the generator $\mathcal{A}$ defined by \cref{def A} verifies all the hypotheses of  the Hille-Yosida Theorem. 
	First,  $\mathcal{A}$ is defined upon the set  $\mathcal C^1_b(\R)$, a dense part of $\mathcal C_b(\R)$. 
	
	Next, the operator $\mathcal{A}$ verifies the maximum principle.  This is straightforward from \cref{def A}, but  maybe even clearer on the following equivalent expression of $\mathcal A$ (see also \cref{pde w}):
	\begin{equation}
	\mathcal A \psi (z)= \beta \ds \int_{\R} K(h) \frac{F(z + \sigma h)}{F(z)}\Big( \psi(z+\sigma h)-\psi
	(z) \Big)dh + c\, \p_z \psi (z), \quad z \in \R.
	\end{equation}
	
	Finally, one needs to check that there exists $\theta \in \R$  such that for any $g$ in a dense subset of $\mathcal C_b(\R)$, the equation
	\begin{align}\label{eqint link}
	\theta \psi - \mathcal A \psi = g
	\end{align}
	admits a solution $\psi  \in \mathcal{C}_b^1(\R)$.
	This is a direct consequence of a similar statement for $\mathcal{L}$, for which a justification of this standard result can be found in \cite{bansayecloezgabriel_ergodic}. Remember that  $\mathcal L$ is introduced in \cref{def L frac}, in particular, $F$ is a stationary profile relatively to $\mathcal L$,
	\begin{align}\label{F stat}
	\mathcal L(F)=0.
	\end{align}
	In addition, $\mathcal A$ and $\mathcal L$ satisfy
	\begin{align}\label{link AL}
	\mathcal A \psi = \frac{\mathcal L(F \psi)}{F}, \quad \text { for all } \psi \in \mathcal C^1_b(\R).
	\end{align}
	To see this, use \cref{def L} to write
	\begin{align}\label{eqint AL}
	\frac{1}{F} \mathcal{L}(F \psi) =  \frac{\beta}{F} \Big( K_\sigma* (F \psi)-F \psi \Big)  + \frac 1F c \,\p_z (F \psi)  - \frac 1F \Big(\mu - \bar \mu \Big) \psi  F.
	\end{align}
	Multiplying each side of the equality \cref{F stat} by $\psi/F$, one obtains:
	\begin{align*}
	\frac {\psi} F c\, \p_z (F)  - \frac 1F \Big(\mu - \bar \mu \Big) \psi F -\frac 1F \beta \psi F = - \frac { \psi } F \beta K_\sigma *F .
	\end{align*}
	Plugging this into \cref{eqint AL},
	\begin{align*}
	\frac{1}{F} \mathcal{L} (F \psi) =  \frac{\beta}{F} \Big( K_\sigma * (F \psi)- \psi  K_\sigma *F \Big)  +  c\, \p_z\psi = \mathcal A \psi .
	\end{align*}
	This proves \cref{link AL} and \cref{eqint link} as a result.
	Finally, the Hille-Yosida Theorem  \cite[\textcolor{o}{Chapter 4, Theorem 2.2}]{ethierkurtz}, proves the existence of a strongly continuous semigroup $(\mathcal M_s,s \geq 0)$ satisfying for every $\psi \in \mathcal C_b(\R)$ and $s>0$:
	\begin{align*}
	\frac{d}{ds} \mathcal M_s \psi  = \mathcal M_s \mathcal A \psi, \quad \text{ and } \ \mathcal M_0  \psi =\psi.
	\end{align*}
(2)	To conclude, we need to prove the relationship \cref{def M}. For all $(t,z) \in \R_+ \times \R$, $0\leq s \leq t$ and $k\in \N$, we  compute:
	\begin{align}\label{eqint link2}
	\frac{d}{ds} \mathcal{M}_{t-s} \left( \frac{\upsilon^k(s,\cdot)}{F}\right)(z) = - \mathcal{M}_{t-s} \left[ \mathcal A \left(  \frac{\upsilon^k(s,\cdot)}{F} \right) \right](z) + \mathcal{M}_{t-s} \left[ \frac{ \p_s \upsilon^k(s,\cdot)}{F}  \right](z).
	\end{align}
	We plug the identity \cref{link AL} for the first term, and \cref{def L frac} for the second. This yields
	\begin{align*}
	\frac{d}{ds} \mathcal{M}_{t-s} \left( \frac{\upsilon^k(s,\cdot)}{F}\right)(z) &  = - \mathcal{M}_{t-s} \left[ \frac { \mathcal L \left( \upsilon^k(s,\cdot) \right) }{F}\right](z) + \mathcal{M}_{t-s} \left[ \frac{ \mathcal L(  \upsilon^k(s,\cdot) ) }{F}  \right](z),\\
	& = 0.
	\end{align*}
	Since $\mathcal M_0 = Id$, by evaluating in $s=t$ we have shown that for all $k \in \N$ and $(t,z) \in \R_+ \times \R$,
	\begin{align}
	\frac{\upsilon^k(t,z)}{F(z)} = \mathcal{M}_{t-s}\left(\frac{\upsilon^k(s,\cdot)}{F}\right)(z), \  \, 0 \leq s \leq t.
	\end{align}
	The identity \cref{def M} is now established.
	\qed
	
	As a matter of fact, we  can also straightforwardly define the Markov process $(Y_s, s\geq 0)$ associated to $(\mathcal M_s,s\geq 0)$, see for instance \cite[\textcolor{o}{Chapter 4, Theorem 2.7}]{ethierkurtz}, and  the dual relationship  \cref{dual process} is a straightforward consequence of \cref{def M}.
	
	\section{Long time behavior of lineages}\label{sec:long_time_lineage}
	We now study of the long time behavior of the ancestral process $Y_s$. We can observe from the duality relationship \cref{dual process} that the time of ancestors is \textit{``backwards''} compared to  the time $t$ of the equilibrium $F$ (and the fractions, per \cref{ass subpop}). The regime $s \to + \infty$ corresponds to the study of the most ancient ancestor.  
	\subsection{Long time behavior of neutral fractions}
	\subsubsection*{Proof of \cref{long time asymp}}$ $\\
	Let  $\psi \in \mathcal{C}^1_b(\R)$, then  $\upsilon_0 := F \psi $ is also a function of $\mathcal{C}^1_b(\R)$. We then consider a neutral fraction $\upsilon$   with initial condition $\upsilon_0$. With \cref{eq:f_sub,ass subpop}, this means
	\begin{equation}\label{sys frac int}
	\left\{ \begin{array}{l}
	\ds \p_t \upsilon (t,z) = \mathcal{L}(\upsilon(t,\cdot))(z) \text{ for } \, t>0 , \, z\in\R\\
	\upsilon(0,z) = \upsilon_0(z) , \ z\in\R,
	\end{array}\right.
	\end{equation}
	and $\mathcal L$ is defined by \cref{def L}. The property \cref{dual process} applied at  $s=t$ yields 
	\begin{align}\label{eqint longtime}
	\frac{\upsilon  (t,z)}{F(z)} = \mathbb{E}_z\left[ \frac{ \upsilon_0 (Y_t)}{F(Y_t)}\right]=\mathbb{E}_z\left[  \psi (Y_t) \right].
	\end{align}
	From \cite[\textcolor{o}{Theorem 2.1}]{gabrielcloez2019}, where the long time behavior of solutions of~\cref{sys frac int} is studied, we obtain 
	\begin{align}\label{asymp upsilon} 
	\upsilon(t,\cdot) \xrightarrow[t \to \infty]{L^\infty_{loc}} p[\upsilon_0] F,  \ \hbox{ with } \ p[\upsilon_0]:= \dfrac{\ds\int_\R \upsilon_0(z)\varphi(z)dz}{\ds\int_\R F(z')\varphi(z')dz'},
	\end{align}
	where $\phi$ is defined as the  solution of  \cref{def dual L}.
	As a result, we have
	\begin{align*}
	p[\upsilon_0] = \dfrac{\ds\int_\R F(z)  \psi (z)\varphi(z)dz}{\ds\int_\R F(z')\varphi(z')dz'} := \mathbb E_z \left[  \psi (Y_\infty)\right]
	\end{align*}
	Therefore, with \cref{eqint longtime}, we have shown that for any $\psi \in \mathcal C_b^1(\R)$,
	\begin{align*}
	\mathbb{E}_z[\psi(Y_t)] \xrightarrow[t\to \infty]{} \mathbb E_z \left[  \psi (Y_\infty)\right] . \end{align*}
	\qed
	
	This proof  hinges  on the recent results of \cite{gabrielcloez2019}, and we believe it can be deduced  as well from the general semigroup analysis of growth fragmentation equations presented in  \cite{mischlerscher}.  In any case, the linearity of the operator $\mathcal{L}$ is crucial to the argument. Indeed, equation \cref{asymp upsilon} states the convergence towards the projection on the dual eigenspace generated by $\varphi$ solution of \cref{dual process} as well as $\mathcal L ^\ast (\phi)= 0$.

	Moreover, the asymptotic proportion of the fractions coincides with the 
	heuristics proposed in \Cref{sec heur pde}. Indeed, choosing $\upsilon_0^y (z)= \delta(z-y)F(z)$ in \cref{sys frac int} (as in \cref{eq:f_sub_dirac}), and applying our convergence theorem~\cref{asymp upsilon},  we recover '$p[y]$', which corresponds exactly to the asymptotic density stated in \cref{long time asymp}.

	\subsubsection*{Proof of \Cref{cor explicit asymp}}\label{sec cor longtime}$ $\\
	The key is to find out that there exists an explicit link between the solution of the dual problem $\varphi$ and the original one \cref{eq:equilibrium c}. Thanks to the symmetry of $\mu$, the function  $ z  \mapsto F(-z)$ solves the dual problem~\cref{def dual L}.  Therefore, the function $F\varphi$ is even, and $Y_\infty$ admits an even density given in \Cref{long time asymp}. \Cref{cor explicit asymp} immediately follows. \qed
	

	\subsection{Numerical simulations : dynamics of the ancestral process: }\label{sec num}$ $\\
	In this section we present some numerical simulations illustrating the previous results. 
	\begin{figure}
		\centering
		\subfigure[Initial Time]{	\includegraphics[scale=0.54]{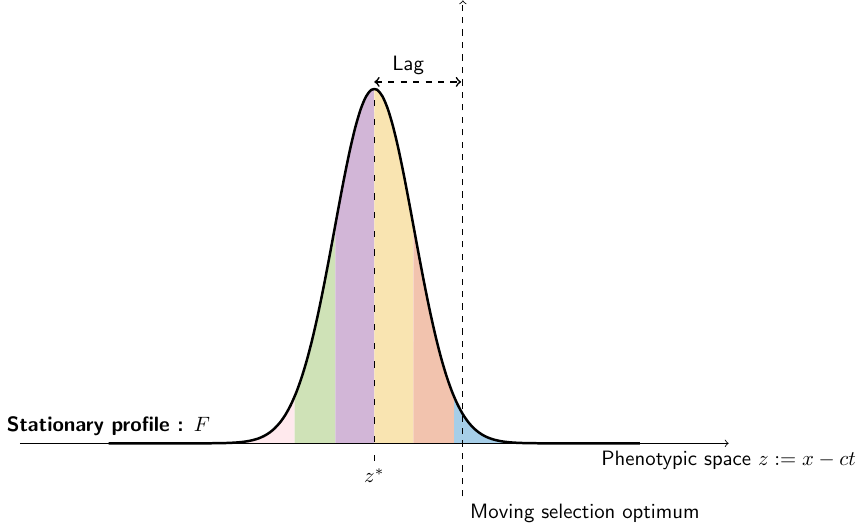} \label{fig initial frac} }
		\subfigure[After a long time]{		\includegraphics[scale=0.54]{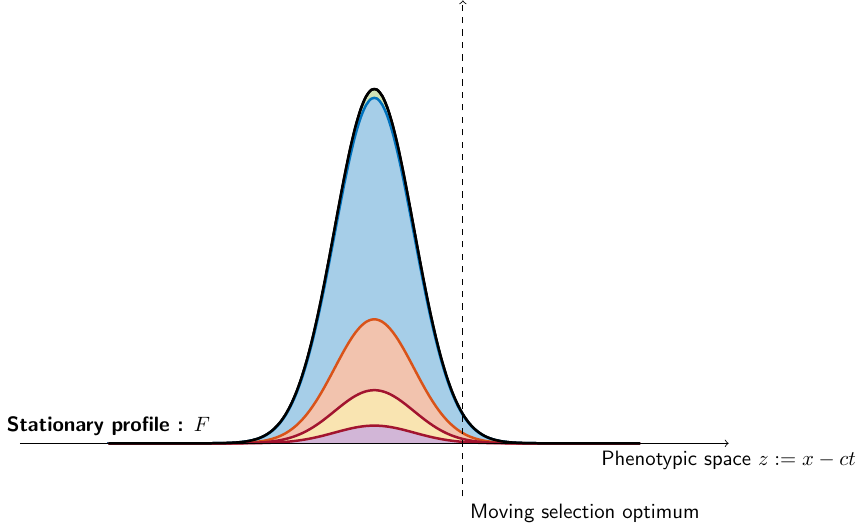}  \label{fig long frac} }

		\caption{Dynamics of the neutral fractions inside the population in the moving frame at speed $c$. The black curve corresponds to the constant profile $F$ defined by~\cref{eq:TW_c}. Each colour corresponds to a different fraction which only differ by their initial location in the moving frame.}\label{fig neutral frac}
	\end{figure}
	\Cref{fig neutral frac} shows  the initial and 'final' distribution of fractions solving \cref{eq:f_sub}. Here, each fraction corresponds to a given interval of traits among the population, and each color corresponds to a neutral fraction. Thus, the cumulative density of the neutral fractions stays equal to the entire population represented by the stationary density profile $F$ moving at constant speed  ( per \cref{ass subpop}). The dominant fraction inside the population after a long time is the 'blue' fraction, which corresponds to the few individuals optimally adapted initially. Thus, we see that the blue fraction contributes the most to the adaptation of the population to the changing environment.  
	
\begin{figure}
		\centering
		\subfigure[Ancestral process $Y_s$ for different times $s$]{	\includegraphics[scale=0.4]{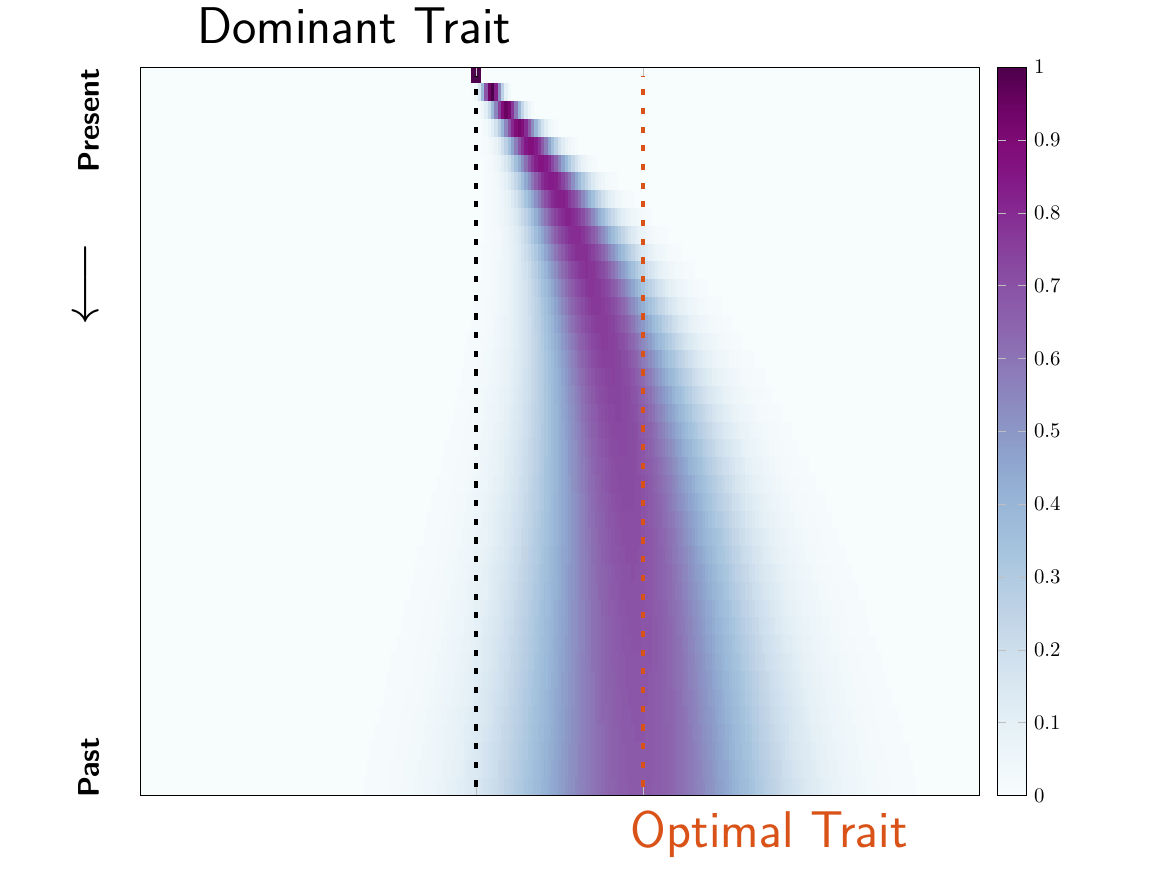}  \label{fig ancestral}}
		\subfigure[In solid blue: density of $Y_\infty$, and in dashed  black: the equilibrium $F$.]{	\includegraphics[scale=0.35]{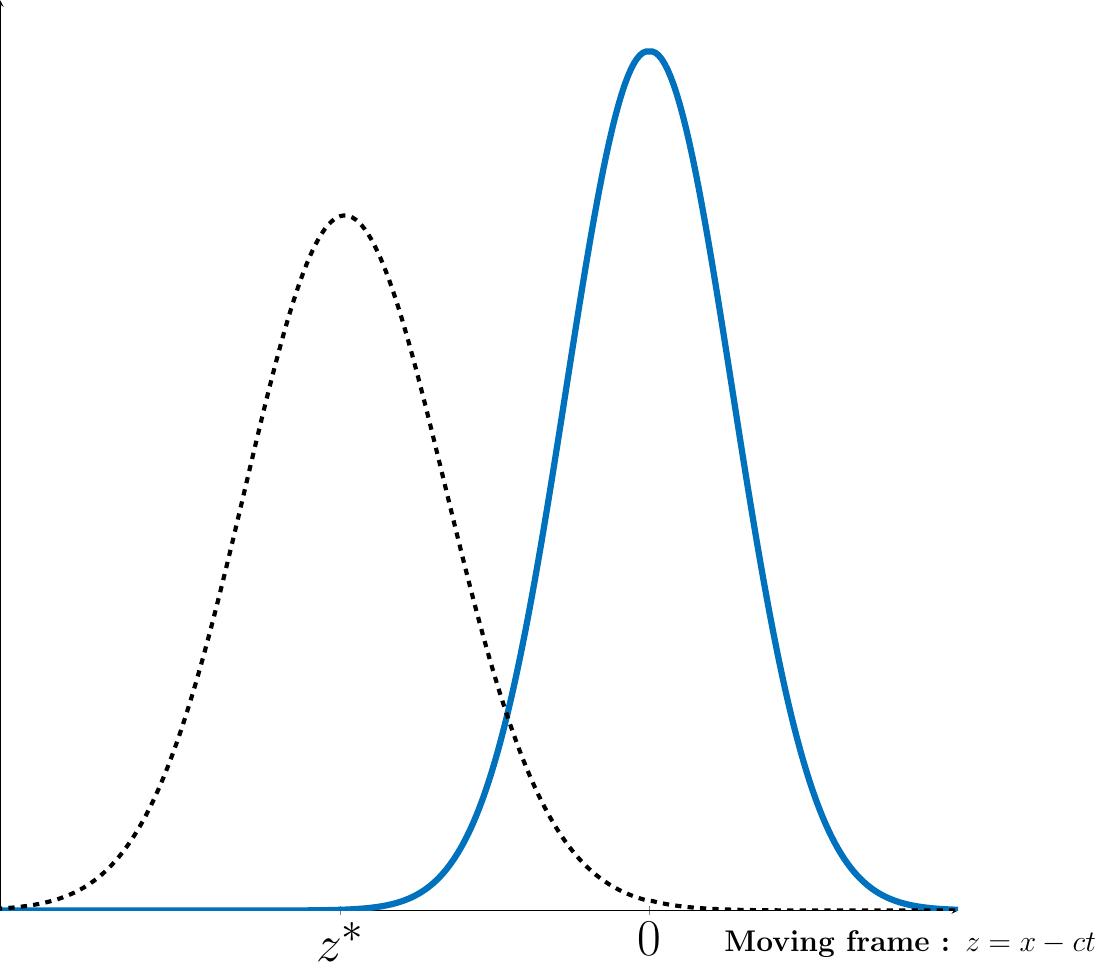} \label{fig long ancestral}}

		\caption{Ancestral Process distribution over time  in panel (a) and asymptotic ancestral distribution $Y_\infty$ (blue curve) compared with the population density $F$ (dotted black curve) in panel (b). In panel (a) the black dotted line correponds to the dominant trait $z^*$ inside the population and the red dotted curve represents the optimal trait $z=0$.}\label{fig fractime}
	\end{figure}

	We now turn to the numerical simulations of the ancestral process, introduced in \Cref{theo edp lineages}. In  \Cref{fig ancestral},  we represent the dynamics of the distribution of the ancestral process $Y_s$ starting at $s=0$ from a single trait $z=z^*$ corresponding to the dominant trait in the population, that is $\max(F) = F(z^*)$.
	We can observe that the ancestral process distribution starts from a Dirac mass at $z=z^*$ when $s=0$ (top), and then it flattens over $z$ and its mean gradually shifts towards the optimal trait $z=0$ as $s$ increases. Eventually, it reaches a stationary distribution which corresponds to the explicit expression of~\cref{long time asymp}
	(see blue curve in \Cref{fig long ancestral}). This asymptotic density represents the  proportion of ancestors of phenotype $y$ in the population, asymptotically as $s \to \infty$. We recover that the asymptotic density is even with a maximum at $z=0$, the optimal trait. Although the density at the optimal trait is very low (see dashed curve in \Cref{fig long ancestral}), most common ancestors, which are ancestors when $s\to\infty$, have a trait close to this optimal trait (see blue curve in \Cref{fig long ancestral}).

	\section{The diffusive approximation}\label{sec diffusive approx}$ $\\
	When the mutational variance $\sigma^2$ is small, the convolution operator $\mathcal{B}$ defined in \cref{eq:asexual B} can be approximated by:
	\begin{equation*}
	\mathcal{B}(f)(z) \approx  f(z) + \dfrac{\sigma^2}{2}\partial^2_z f(z).
	\end{equation*}
	In the diffusive approximation regime, the profile $F_d$ satisfies the following equation:
	\begin{equation}\label{eq:TW_c_diff}
	0 = \beta\dfrac{\sigma^2}{2}\partial^2_z F_d(z) + c\partial_z F_d(z) + \left(\mu(z)-\gamma \int_\R F_d(z')dz'\right)F_d(z).
	\end{equation}
	Therefore,  the neutral fractions  operator $\mathcal{L}$ defined by \cref{def L} becomes:
	\begin{equation}\label{def L diff}
	\mathcal{L}_d(\upsilon) = \beta\dfrac{\sigma^2}{2}\partial^2_z \upsilon + c\partial_z \upsilon - \Big(\mu -\overline{\mu}  \Big)  \upsilon 
	\ \hbox{ with } \ \bar\mu = \gamma \int_\R F_d(z)dz = \int_\R\mu(z) \dfrac{F_d(z)}{\int_\R F_d(z')dz'}dz.
	\end{equation}
	With this new model, we can define the ancestral lineage under the diffusive approximation and we can state the following properties. 
	\begin{prop}[Ancestral process under the diffusive approximation]\label{prop diff}$ $\\
		(a) The \textbf{diffusive ancestral process} $(Y_s^d, s\geq 0)$ associated to  \cref{def L diff} admits the following generator:
		\begin{align}\label{def A diff}
		\mathcal A_d \psi = \beta \frac{\sigma^2}{2} \p_y^2 \psi  + \left( \beta \frac{\sigma^2}{2} \frac{\p_y F_d} {F_d}   + c \right)  \partial_y \psi.
		\end{align}
		(b) \textbf{Diffusive common ancestor distribution}. When $s \to \infty$, the limit process $Y_\infty^d$ admits a density given by
		\begin{align}\label{diff asymp dens}
		\dfrac{\ds \left( F_d(y) e^{cy /(\beta\sigma^2)} \right)^2}{\ds\int_\R \Big(F_d(y')e^{cy'/(\beta\sigma^2)}\Big)^2dy'}.
		\end{align}
		(c)  In addition, the above density admits a local maximum at $y=0$ for $\sigma $ small enough.
	\end{prop}
	In this case, equation~\cref{def A diff} implies that the ancestral process $(Y_s^d, s \geq 0) $  solves the following SDE:
	\begin{align} \label{SDE_Y_diffusive}
	d Y_s = \left( \beta \frac{\sigma^2}{2} \frac{\partial_y F_d(Y_s)}{F_d(Y_s)} + c \right) ds + \sqrt{\beta \sigma^2} dB_s,
	\end{align}
	where $ (B_s, s \geq 0) $ is standard Brownian motion. A similar statement is made in \cite[\textcolor{o}{Theorem 1.1}]{companion} in the case of quadratic selection ($F_d$ is then explicit, see below).
	
	\subsubsection*{Proof of \Cref{prop diff}.}
	
	(a)~Formula \cref{def A diff} is obtained by plugging expression \cref{def L diff} in  \cref{link AL}, yielding
	\begin{align*}
	\mathcal A_d \psi  =  \frac 1{F_d} \mathcal L_d (F_d \psi) = \beta \frac{\sigma^2}{2} \left( \frac 1{F_d} \p_y^2 (F_d \psi ) - \frac 1{F_d} \psi \, \p_y^2 F_d \right) + c \partial_y \psi. 
	\end{align*}
	Then, after simplifications using \cref{eq:TW_c_diff}, we recover \cref{def A diff}. 
	
	\
	
	(b)~To establish \cref{diff asymp dens},  we introduce the dual problem, similarly to  \cref{def dual L}:
	\[ \mathcal{L}_d^\ast(\phi) = \beta\dfrac{\sigma^2}{2}\partial^2_y \phi - c\partial_y \phi  - \Big(\mu -\overline{\mu}  \Big)  \phi.
	\] 
	In order to prove  the convergence of $Y_s$, we must use similar result to \cref{asymp upsilon}, with $\mathcal L _d$. This is a classical statement, see for instance \cite{garnier2012inside} and references therein (the convergence now holds in $L^2$). 
	Moreover the density in \cref{diff asymp dens} is slightly more explicit than in \cref{long time asymp}, because there exists a special relationship between the primal and dual problem. Indeed, computations show that  if $\mathcal{L}_d^\ast(\phi_d) = 0$, then
	\begin{align*}
	\varphi_d(y) = F_d(y)e^{2cy/(\beta\sigma^2)}.
	\end{align*}
	This explains \cref{diff asymp dens}, using the formula of  \cref{long time asymp}. 
	
	\
	
	(c)~The function defined for each $z$ by $\tilde F(y) := F_d \, e^{cy/(\beta\sigma^2)}$ is even, since $\mu$ is an even function by hypothesis and $\tilde F$ is the  solution of: 
	\begin{equation}\label{eq tilde F}
	\beta\dfrac{\sigma^2}{2}\partial^2_y \tilde F - \Big(\mu -\overline{\mu} + \dfrac{c^2}{2\beta\sigma^2}  \Big)  \tilde F= 0.
	\end{equation}
	Therefore $y \mapsto \Big( F_d(y) e^{cy /(\beta\sigma^2)}\Big)^2$ is also an even function. As a consequence $Y_\infty$ is symmetric, by \cref{diff asymp dens}. In particular, $\tilde F$ admits a (local) extrema at $z=0$. To obtain more information  on the local shape of the density, we must investigate the sign of $\p_y^2 \tilde F(0)$, since $\tilde F$ is the numerator in \cref{diff asymp dens}. First,  we find that $\mu(0)-\bar \mu <0$. Back to \cref{eq tilde F}, if $\sigma$ is sufficiently small, we get  $\p_y^2 \tilde F(0)>0$, and therefore  $Y_\infty$ admits a local maximum at $0$, as claimed in \Cref{prop diff}. 
	\qed
	
	\subsection*{Quadratic selection. }
	Here, we tackle the special case where $\mu(z) = z^2/2$,  which verifies our general hypotheses in \cref{hyp growth m}. In this case, the profile $F_g$ is known to be Gaussian, which allows us to illustrate our previous qualitative comments~\cite{BurLyn95,kopp2014rapid}. More precisely, the profile $F_g$ solves 
	\begin{equation}\label{eq:diff limit}
	\lambda_g  F_g(z) - c \partial_z F_g(z) + \frac{z^2}2  F_g(z)=  \beta F_g(z) + \dfrac{ \beta \sigma^2}{2}\partial_z^2 F_g(z)  , \text{ for}  z\in \R,
	\end{equation}
	whose solution is given by
	\begin{equation} \label{F_quad_sel}
	F_g(z) = \dfrac{\lambda_g}{\sqrt{2\pi\sigma\sqrt \beta }}\exp\left ( -\frac{1 }{2\sigma\sqrt \beta}  \left ( z + \frac c{\sigma \sqrt \beta}  \right )^2 \right )\, , \quad \lambda_g = \beta - \dfrac{c^2}{2 \beta \sigma^2 } - \dfrac{\sigma\sqrt \beta}{2} \, .
	\end{equation}
	Up to a constant, $F_g$ is indeed a Gaussian distribution centered around an optimum proportional to $c$,  lagging behind the optimal trait, and with a variance proportional to $\sigma$ (instead of $\sigma^2$ for $K_\sigma$). In the eigenvalue $\lambda_g$, we recognize the \emph{lag load}: $c^2/2 \beta \sigma^2$, that is the cost to keep pace with the changing environment, and the \emph{mutation load}: $\sigma \sqrt \beta /2$ \cite{BurLyn95}. In addition, the speed of change $c$ must be small enough for the population to persist ($\lambda>0$ if and only if $c < \sigma\beta\sqrt{2}\sqrt{1-\sigma/(2\sqrt{\beta})}$) \cite{Bur00}.
	In particular, $c$ must be of the order of $\sigma$, just as in the regime of \Cref{regime adaptative dynamics}.
	
	Moreover, from~\cref{diff asymp dens}, we get the following explicit formula for the asymptotic  ancestral distribution $Y_\infty^g$:
	\[
	\frac{ (F_g\phi_g)(y)} {\ds \int_\R (F_g\phi_g)(y')\, dy'}= \dfrac{1}{\sqrt{\pi\sigma\sqrt \beta }}\exp\left( -\frac{y^2 }{\sigma\sqrt \beta} \right). 
	\]
	This is a Gaussian distribution centered at $y=0$ and with variance $\sigma\sqrt{\beta}/2$. The variance of the asymptotic ancestral distribution is half that of the population and coincides with the mutation load.
	
	In addition, under the diffusive approximation, we can characterize the entire ancestral process $ (Y_s, s \geq 0) $ by combining the expression of $F_g$ stated in~\cref{F_quad_sel} and the SDE equation solved by $Y_s$.
	The ancestral process is an Ornstein-Uhlenbeck process. This characterization can also be found in \cite[\textcolor{o}{Theorem 1.1}]{companion} (with $\beta=1$, \cref{F_quad_sel} corresponds to \cite[\textcolor{o}{Proposition 2.3}]{companion}). 
	Therefore, for all $s \geq 0$, $Y_s$ follows a Gaussian distribution with mean and variance given by
	\begin{align}\label{mean var quadex}
	\forall s\geq 0, \quad \mathbb{E}_z(Y_s) = ze^{-\sigma\sqrt{\beta}s} \ \hbox{ and  } \ 
	\Var_z(Y_s) = \dfrac{\sigma\sqrt{\beta}}{2}\left(1-e^{-2\sigma\sqrt{\beta}s}\right).
	\end{align}
	The variance does not depend on the reference point $z$, while the mean of the ancestral distribution converges to $0$ exponentially fast, at a rate $\sigma\sqrt{\beta}$.  
	Building on this, we  conjecture that for general symmetric selection function $\mu$, the  mean $\mathbb{E}_z (Y_s)$ should converge exponentially fast to $0$ at a rate $\sigma^2\beta/\Var(F_d)$, where $\Var(F_d)$ corresponds to the variance of the profile $F_d$ defined by \cref{eq:TW_c_diff}, we refer to \Cref{sec approx app} for further details.

	

	\section{Small mutations regime}\label{sec adap}
	In this section, we tackle the asymptotic regime where mutations have very small effects (the trait of a mutant offspring is very close to that of its parent). First, we explain how  we can deduce the heuristic formula \cref{edo Gamma} for the typical lineage from the following Hamilton-Jacobi equation obtained as the limit when $\sigma \to 0$ of the equation satisfied by $F$:
	\begin{align}\label{spec Ham}
	\lambda  + c' \partial_{z} U(z)  + \mu(z) = \beta + \beta H(\p_z U(z)) .
	\end{align}
	where $U=\lim_{\sigma\to0}{\sigma \ln(F)}$ and $H$ is the Hamiltonian defined in~\Cref{regime adaptative dynamics}.
	The rigorous convergence when $\sigma \to 0$, 
	in the sense of viscosity solutions, is established in \cite{barles2009,lorz2011dirac}.
	In \Cref{sec HJ derivation}, we explain how to  derive (formally) the Hamilton-Jacobi equation from \cref{eq:equilibrium c} when $\sigma \to 0$.  
	
	\subsection{Hamilton-Jacobi equation and the typical lineage}\label{sec HJ}$ $\\
	Before proving \Cref{regime adaptative dynamics}, we start with some heuristic arguments linking the Hamilton-Jacobi equation~\cref{spec Ham} and the typical lineage using a dual point of view.  
	Let us introduce the Lagrangian function $L$ associated to $K $, corresponding to the Legendre transform of the Hamiltonian $H$:
	\begin{align}\label{def Lag}
	L(v) := \max_{p \in \R} \Big(  pv- H(p) \Big).
	\end{align}
	Using this function, we can write the solutions of \cref{spec Ham} with this variational representation:
	\begin{align}\label{asex var}
	U(z) = \inf_{\gamma \text{ s.t. } \gamma(0)=z} \int_{0}^{+ \infty} \left[ \beta L\left(\frac{-\dot{\gamma}(s) + c'}{\beta}\right) - \beta + \mu(\gamma(s)) + \lambda \right] ds.
	\end{align}
	We refer to \cite{barles2006} for the origin of this formula, which stems from the Weak-KAM theory. 
	 The infimum is taken over all functions $\gamma \in \mathcal C^1(\R_+)$ that reach the phenotype $z$ at time $0$. 
These functions $\gamma$ are phenotypic paths, among them  any optimal trajectory  minimizes a cost. 
	The expression of this functional shows the combined cost of mutations through $L$ (at speed $\dot \gamma +c'$) and  selection through $\mu$. 
	The birth rate $\beta$ plays a role opposite to natural selection, while $\lambda$ is the term that balances the expression, just as in \cref{spec Ham}.
	
	\textbf{From a genealogical point of view, the phenotypic path can be seen as an ancestral lineage. As a result, the special lineage that minimizes the cost should correspond to a 'typical'  ancestral lineage of the individuals of trait $z$ at time $0$.} 
	
	
For the time being, let  $\Gamma$ be such a minimizing trajectory.  Then, as a byproduct of the Weak-KAM theory, on can show that  the optimal trajectory $\Gamma$ of the variational problem \cref{asex var} is the solution of  the following ordinary differential equation:
	\begin{align}\label{edo Gamma_bis}
	\dot \Gamma(s)& = c' - \beta \, \p_p H \Big( \p_z U \big( \Gamma(s) \big) \Big) , \quad s>0\\
	\nonumber \Gamma(0) & =z .
	\end{align}
	This result comes from the Hamiltonian/Lagrangian structure of \cref{asex var} and more precisely from the study of the characteristics of this Hamilton Jacobi equation (see for instance \cite{hairer2006}). 
	This characterization of the optimal trajectory will help us gain qualitative insights on the typical lineage and will improve our numerical computation of $\Gamma$~\cite{hairer2006}. 
	
	
The ODE \cref{edo Gamma_bis} coincides with the limit equation~\cref{edo Gamma} of the ancestral process in \Cref{regime adaptative dynamics}. This is the key point of this part, and as a result,  we can state qualitative properties of the ancestral process studying $\Gamma$. 
	If $\Gamma$ satisfies~\cref{edo Gamma_bis} then we have
	\begin{align}\label{U0 Gamma}
	U(z) = \int_{0}^{+\infty}  \left[ \beta  L\left(\frac{-\dot{\Gamma}(s) + c'}{\beta}\right) - \beta + \mu(\Gamma(s)) + \lambda \right] ds.
	\end{align}
	Formally, one should expect $\dot \Gamma(s)$ to converge to  $0$ when $s \to +\infty$. 
	Otherwise, $\Gamma$ and thus the mortality  rate $ \mu (\Gamma)$ would become arbitrarily large ( see \cref{hyp growth m}). Thus, the trajectory $\Gamma$ would not be a  minimizer, as $U(z)$ could take the value $+\infty$. As a result, $\Gamma$ converges when $s \to +\infty$ and its limit is necessary $0$  to minimize the selection function $\mu$:
	\begin{align}\label{lim gamma}
	\Gamma(s) \xrightarrow[s \to +\infty]{} 0.
	\end{align}
	The typical lineage traces back to an ancestor with the optimal trait $0$ as we already observed in \Cref{cor explicit asymp};  however, here it is in a much stronger sense.
	We prove the limit \cref{lim gamma} in \Cref{sec lim Gamma},  using convexity methods and qualitative properties extracted from the Hamilton-Jacobi equation beyond lineages (see \Cref{sec Lag pov}). 


	\subsection{Proof of \cref{regime adaptative dynamics}}$ $\\
	We aim to prove the convergence of the semigroup $\mathcal{M}^\sigma$ as $\sigma\to0$. From classical results on semigroups, we only need to prove the convergence of its  generator $\mathcal{A}^\sigma$~\cite[\textcolor{o}{Theorem 17.25}]{kallenbergfoundations}.
	We  prove the convergence of the generator on $\mathcal{C}^2_b(\R_+ \times \R )$ which is actually a core for the generator $\mathcal{A}^\sigma$~ \cite{bansayecloezgabriel_ergodic,gabrielcloez2019}. 
	Now, let $\psi \in \mathcal{C}^2_b(\R_+ \times \R )$ and $\mathcal A^\sigma$ be the generator corresponding to the semigroup $ \mathcal M^\sigma $. Then, we have
	\begin{equation*}
	\mathcal A^\sigma \psi(y) = \frac \beta \sigma \ds \int_{\R} K (h) \frac{F(y + \sigma h)}{F(y)}\Big(  \psi(y+\sigma h)-\psi (y) \Big)dh + \frac c \sigma  \p_y \psi(y), \quad y \in \R.
	\end{equation*}
	Let $U_\sigma$ be defined by  \begin{align}\label{hopf asex bis}
	U_\sigma(z) = - \sigma  \log( F_\sigma(z)) \ \hbox{  for all } \ z \in \R .
	\end{align}
	Then the generator $ A^\sigma$ is written as
	\begin{equation*}
	\mathcal A^\sigma \psi(y) = \frac \beta \sigma  \ds \int_{\R} K (h)  \exp \left( -\frac{ U_\sigma(y + \sigma h) - U_\sigma(y)}{\sigma} \right) \Big(  \psi(y+\sigma h)-\psi(y) \Big)dh + \frac c \sigma \p_y \psi(y).
	\end{equation*}
	Using a Taylor expansion, for all $y \in  \R$, there exists $\tilde y$ such that $ |y-\tilde y| \leq \sigma h$ and
	\begin{align*}
	\psi(y+\sigma h)-\psi (y) = \sigma h \p_y  \psi(y) + \frac{\sigma^2 h^2 }{2} \p_y^2 \psi(\tilde y).
	\end{align*}
	Plugging this expression into the definition of the generator $\mathcal{A}^\sigma$, we end up with
	\begin{multline}\label{eqint lim}
	\mathcal A^\sigma \psi(y) = \beta  \ds \int_{\R} K(h)  \exp \left( -\frac{ U_\sigma(y + \sigma h) - U_\sigma(y)}{\sigma} \right) h \, \p_y  \psi(y)  \, dh   \\ + \sigma \beta \ds \int_{\R} K(h)  \exp \left( -\frac{ U_\sigma(y + \sigma h) - U_\sigma(y)}{\sigma} \right) \frac{ h^2 }{2} \p_y^2 \psi(\tilde y) \, dh  +  \frac c \sigma \p_y \psi(y).
	\end{multline}
	Thanks to \cite{barles2009,lorz2011dirac}, we know that $U_\sigma$ converges locally uniformly towards $U$ as in \cref{limit U}. 
	Moreover $U$ solves the problem \cref{spec Ham}, in the sense of viscosity solutions, and we have the following Lipschitz uniform bound:
	\begin{align} \label{Lipschitz}
	\norm{\,U' _\sigma \, }_\infty <\eta,
	\end{align}
	for all $\sigma > 0$ small enough, where $ \eta $ is such that
	\begin{align*}
	\int_\R K(h) e^{\eta |h|} dh < +\infty,
	\end{align*}
	As a result, we show that the first intergal term of \cref{eqint lim} converges and the second integral vanishes, pointwise, for any $\psi \in \mathcal{C}^2_b (\R_+ \times \R)$.
	Using the rescaled speed $c' = c /\sigma$, we obtain
	\begin{align}\label{def Azero}
	\mathcal A^\sigma \psi(y) \xrightarrow[\sigma \to 0]{} - \beta  \ds \int_{\R} K(h)  \exp \Big(h \,  U' (y) \Big)h \, \p_y  \psi(y) \,   dh  + c' \p_y \psi(y):= \mathcal{A}^0 \psi(y).
	\end{align}
	Then using Theorem 17.25 in \cite{kallenbergfoundations}, we conclude that the semigroup $\mathcal M _s^\sigma$ associated to $\mathcal{A}^\sigma$ converges to a semigroup $\mathcal{M}_s^0$ associated to the asymptotic operator $\mathcal{A}^0$. 
	
	To conclude the proof of \cref{regime adaptative dynamics}, we need to characterize  $\mathcal M^0$. 
	Let $\mathcal V$ be defined by
	\begin{equation*}
	\mathcal V (z) = c' - \beta \p_p H( U' ( z)  ),  \quad  \text{ for all } z \in \R,
	\end{equation*}
	and $\Gamma$ be its corresponding integral flow: 
	\begin{equation}\label{gamma flow}
	\left\{ \begin{array}{ll}
	\partial_s \Gamma (t,s,z)=\mathcal V (\Gamma(t,s,z)), \\
	\Gamma(t,t,z)  = z.
	\end{array}\right.
	\end{equation}
	The expression for the derivative of the Hamiltonian as defined in \cref{limit U} is:
	\begin{align*}
	\p_p H(p) =  \int_\R  K(y) \exp (y p ) y \, dy.
	\end{align*}
	
	Now for any test function $\psi$, we define $\theta(s,z) = \mathcal M_s^0 \, \psi (z)$.
	Using the expression of $\mathcal A ^0$, we see that $\theta$  solves  the  following equation:
	\begin{equation*}
	\left\{\begin{array}{ll}
	\p_s  \theta (s,z) =  \mathcal V (z) \p_z \theta(s,z), \quad s>0, \,  z\in \R,\\ 
	\theta(0,z) = \psi(z).
	\end{array}\right.
	\end{equation*}
	Classically, this advection equation with non constant velocity field admits an ``explicit'' solution, by following the characteristics, which corresponds in our case to the  flow $\Gamma$ defined by \cref{gamma flow}:
	\begin{equation}\label{eqint limit3}
	\theta(s,z)= \psi(\Gamma(0,s,z)).
	\end{equation}
	From the definition of $\Gamma_z$ stated in \cref{edo Gamma} and the expression \cref{gamma flow} of $\Gamma$, we observe that  $\Gamma(0,s,z) = \Gamma_z(s)$. From the expression of $\theta$ in \cref{eqint limit3}, we prove the convergence of \cref{conv semigroup} in \cref{regime adaptative dynamics}.
	
	\bigskip
	
	We now turn to the proof of \cref{cvg_probability}. We claim that the family of Markov processes $ (Y^\sigma_s, s \geq 0) $ indexed by the parameter $ \sigma \in [0,1] $ is tight for the Skorokhod topology.
	To prove this, note that
	\begin{align*}
	Y^\sigma_s = Y^\sigma_0 + V^\sigma_s + M^\sigma_s,
	\end{align*} 
	where
	\begin{align*}
	V^\sigma_s = c's + \beta \int_{0}^{s} \int_{\R} K(h) \frac{F(Y^\sigma_r+\sigma h)}{F(Y^\sigma_r)} h\, dh dr
	\end{align*}
	and $ (M^\sigma_s, s \geq 0) $ is a local martingale with predictable variation
	\begin{align*}
	\langle M^\sigma \rangle_s = \sigma \beta \int_{0}^{s} \int_{\R} K(h) \frac{F(Y^\sigma_r+\sigma h)}{F(Y^\sigma_r)} h^2 dh dr.
	\end{align*}
	Using \cref{Lipschitz}, we then see that
	\begin{align*}
	&| V^\sigma_{s'} - V^\sigma_s | \leq \left( c' + \beta \int_{\R} K(h) \exp( \|U'_\sigma \|_\infty\, |h|)\, |h| dh \right) |s'-s|, \\
	&| \langle M^\sigma \rangle_{s} | \leq \sigma s \beta \int_\R K(h) \exp( \| U'_\sigma \|_\infty\, |h|)\, h^2 dh.
	\end{align*}
	Since $ \|  U'_\sigma \|_\infty < \eta $ as in \cref{ass K}, there exists a constant $ C > 0 $, independent of $ \sigma $, such that
	\begin{align*}
	\int_{\R} K(h) \exp(\| U'_\sigma \|_\infty\, |h|)\, |h|\wedge |h|^2 dh \leq C,
	\end{align*}
	where $a\wedge b = \min(a,b)$.
	As a result, for any $\varepsilon > 0$,
	\begin{align*}
	\lim_{\delta \to 0} \limsup_{\sigma \to 0} \mathbb P\left( w_\delta(Y^\sigma) > \varepsilon \right) = 0,
	\end{align*}
	where $w_\delta(Y)$ denotes the modulus of continuity of $Y$, \textit{i.e.}
	\begin{align*}
	w_\delta(Y) = \sup_{\underset{|t-s|\leq \delta}{s, t \in [0,T]}} |Y_s-Y_t|.
	\end{align*}
	Moreover, for each fixed $ s \geq 0 $, $ Y^\sigma_s $ converges in distribution to the deterministic limit $ \Gamma_z(s) $, this convergence holds also for finite-dimensional marginals of $ (Y^\sigma_s, s \geq 0) $.
	This shows that the family of processes $ (Y^\sigma_s, s \geq 0) $ is C-tight \cite{billingsley2013convergence}.
	This yields the convergence in distribution of $ (Y^\sigma_s, s \geq 0) $ in the uniform topology to the deterministic process $ (\Gamma_z(s), s \geq 0) $.
	Finally, since the limit is deterministic, the convergence  holds in probability. Hence \cref{cvg_probability} is proven which concludes the proof of \Cref{regime adaptative dynamics}.
	\qed

	\section{Numerical simulations}$ $\\
	The aim of this section is to compare our ancestral process defined from the deterministic model with the ancestral lineages of classical Individual Based Model (IBM) taken from \cite{champagnat2006unifying,champagnat2007individual}. Moreover, we assess the accuracy of our approximation formula (diffusive approximation and Hamilton Jacobi approximation) with respect to the mean of the ancestral process $Y_s$ and the mean of the lineages of the IBM model.
	
	\subsection{The stochastic model}
	We consider a stochastic IBM model where each individual is characterized by its trait. They reproduce and die at rates that may depend on their traits (see \Cref{app:IBM} for more details) and on the total population size. 
	A logistic competition term keeps the population size finite, and when the strength of competition tends to zero, the population size tends to infinity and the (renormalized) stochastic model converges to the deterministic model \cref{eq:f z}
	 per  ~\cite{champagnat2006unifying} (see \Cref{fig profiles asex} in~\Cref{app:IBM}). 
	Therefore, for large populations, we may expect the lineages of the IBM model and the trajectories of the ancestral process $Y_s$ associated to the deterministic model~\cref{eq:f z} to behave similarly. 
	In order to compare the genealogy of the individuals in the stochastic model and the ancestral process, each individual carries a label which encodes its genealogy, that is the trait of its parents, together with its current trait (see ~\Cref{app:IBM} for more details). 
	With this definition of the genealogy, we compare the dynamics of lineages obtained through the IBM model with the ancestral process distribution for the general model~\cref{eq:f z}, and the lineage trajectory $\Gamma_z$ obtained in the asymptotic regime of small mutations ($\sigma\to0$).
	
\begin{figure}
		\begin{center}
			\includegraphics[scale=0.55]{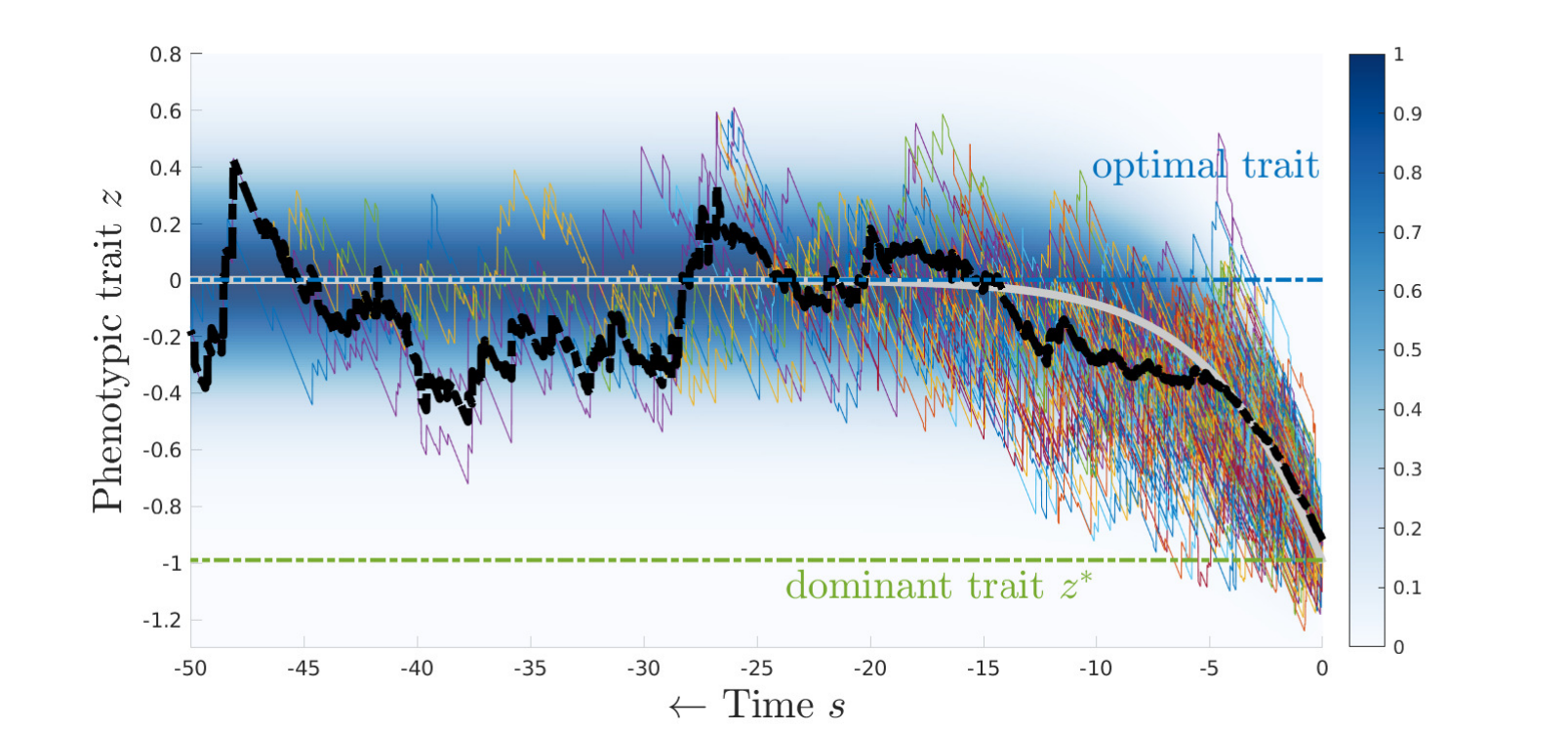}
			\caption{\small Dynamics of the ancestral lineages starting from the dominant trait $z^*$ for the IBM model and the deterministic PDE model. The $y-$axis represents the phenotype in the moving frame. The blue background represents the distribution of the ancestral process $Y_s$ starting from $z^*$ at $s=0$ 
				and the gray curve represents the typical lineage $\Gamma_{z^*}$ defined by~\cref{edo Gamma}.
				The coloured lines starting from $z=z^*$ corresponds to the lineages of the IBM model with $N=2.10^4$ individuals. The black dashed lines represents the mean of the lineages. The blue dash-dotted line represents the optimal trait $z=0$ and the green dash-dotted line the dominant trait $z^*$. For the simulation $\sigma = 0.1$ and $c=0.2.$
			}\label{fig lineage asex}
		\end{center}
	\end{figure}

	From a numerical point of view, we look at the lineages of individuals with the dominant trait of the population $z=z^\ast$ in order to sample initially as many individuals as possible for the IBM model. The population size of the stochastic model is around  20,000 while the number of individuals carrying the dominant trait is around 1,000 (see \Cref{fig lineage asex}). 
	
	\subsection{Ancestral process and stochastic lineages}
	First of all, we have verified that the dynamics of the lineages in both models coincide (see~\Cref{fig lineage asex distrib}  and movie~\ref{app:movie} in \Cref{app:IBM_vs_PDE}). Thus, our ancestral process captures the distribution of the ancestral lineages induced by the stochastic model. 
	In particular, the mean of the ancestral process $\mathbb{E}_z(Y_s)$ and the mean of the stochastic ancestral lineages are close for each replicate (see \Cref{fig lineage asex distrib}). 
	
We compare  the stochastic model with the Hamilton Jacobi approximation obtained in the limit of small mutations ($\sigma\to0$). To solve \cref{edo Gamma} numerically, we coupled this equation with \cref{spec Ham} to obtain a system of ODE. 
	Although $\sigma$ is not $0$ in the stochastic model,  we observe in \Cref{fig lineage asex} that the trajectory $\Gamma_z$ follows the trajectory of the mean of the stochastic ancestral lineages. We further  show in the appendix that the Hamilton Jacobi approximation as well as the  diffusive approximation provide a good estimate of the mean of the ancestral process for both deterministic and stochastic model (see \Cref{fig lineage asex distrib} and \Cref{fig:mean_approx}).
	
\begin{figure}
		\subfigure[$\sigma=0.1$ and $c=0.1$]{\includegraphics[scale=0.5]{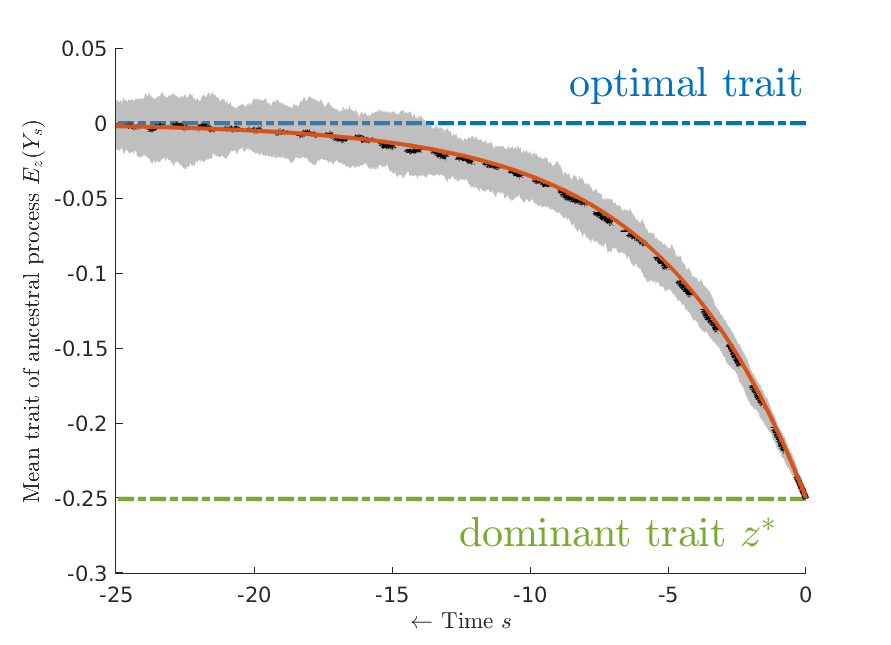} \includegraphics[scale=0.5]{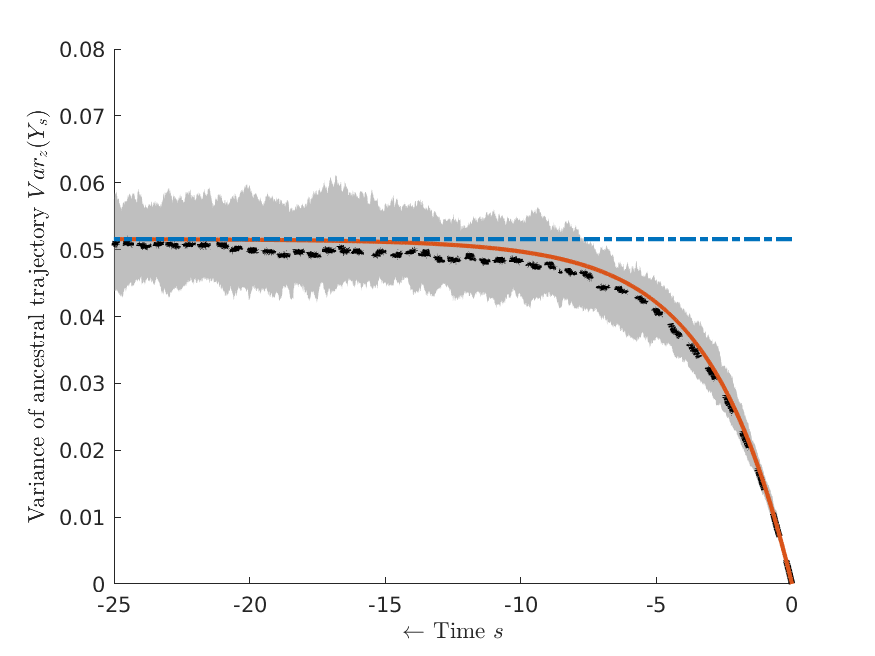}}
		\caption{ \textbf{Mean (left) and variance (right) of the ancestral process $Y_s$ for different set of parameters}: red curves corresponds to the model~\cref{dual process} and the black dashed curves correspond to the IBM model averaged over $50$ replicates and the gray regions corresponds to the $5\%$ and $95\%$ confidence interval (20000 individuals  for each replicate).   
			On the left panel, the horizontal blue line is the optimal trait in the moving frame, $0$. The green line is the dominant trait of the equilibrium $F$, denoted by $z^\ast$. On the right panel, the horizontal cyan line corresponds to the asymptotic variance given by the deterministic model see~\cref{long time asymp}.
		}\label{fig lineage asex mean var bis}
	\end{figure}

	However, going far away in the past in \Cref{fig lineage asex}, the mean of the  stochastic ancestral lineages fluctuates a lot due to the small number of ancestors at this time . In that time regime, it makes more sense to compare averages over a large number of IBM replicates, which we do in  \Cref{fig lineage asex mean var bis}, see also \Cref{fig profiles asex}. 
	As already mentioned the deterministic model is relevant when the size of the population is large enough (see \Cref{fig lineage asex distrib}). However, per \Cref{fig lineage asex mean var bis}, our ancestral process fully captures the average behaviour of the mean over many replicates of the IBM model even if the size of the population is small.

	\section{Discussion}\label{sec disc}

	In the context of adaptation to a changing environment,  
	we proposed  a method to track lineages using a deterministic mathematical model of mutation and selection. More precisely, we are able to define the ancestral process $(Y_s, s \geq 0)$ describing the trajectory of the ancestral traits of an individual sampled uniformly from those with a given trait in the present population. 
	Our results show that every trait is represented among the ancestors, but that ancestral traits are strongly biased towards the fitness optimum (which is shifting linearly).
	Furthermore, combining the asymptotic regime of small mutation ($\sigma \to0$) with the Lagrangian structure of the solution to the Hamilton-Jacobi equation,  we provide a good approximation of the ancestral lineages as the solution to an ODE, see \cref{edo Gamma} in \cref{regime adaptative dynamics}.
	
	Some modeling choices made when writing equation \cref{eq:f z} may seem somewhat arbitrary and limiting the scope of the present study.
	The mathematical analysis performed here can nonetheless be made fairly general.
	In fact, the same analysis can be carried out for any density profile $f(t,\cdot)$ satisfying an equation of the form
	\begin{align*}
	\partial_t f(t,x) = \mathcal{B} f(t,\cdot) (x) + r(f(t,\cdot)) f(t,x),
	\end{align*}
	where $\mathcal{B}^*$ generates a Feller semigroup and which admits a non-negative stationary solution.
	In particular, we could include a dependence of the birth rate on the trait, more general competition terms, etc.
	In fact, we already make fairly weak assumptions on the shape of the function  $\mu$, contrary to, for instance, Fisher's  Geometric Model (FGM) which assumes an explicit quadratic relationship between phenotype and fitness (see~\cite{GilHamMarRoq19,MarRoq16} for related works on the FGM model). 
	In addition, our framework deals with a general form of mutations. 
	Although it can leads to 
	the ``diffusive approximation'' accounting for small mutations~\cite{Kim64} (see~\Cref{sec diffusive approx}),  
	the non local operator $\mathcal{B}$ defined by the mutation kernel $K$ can describe general distributions of mutation effects. 
	We show that the ancestral process truly depends on this mutation kernel (see \cref{pde w}). In particular, we show that, as the effect of mutations vanish ($\sigma\to0$) the lineages retain a trace of the whole mutation kernel through 
	the Hamiltonian in \cref{edo Gamma} (while the diffusive approximation only keeps track of the variance of the mutation kernel through the parameter $\sigma^2$). 
	

	\subsection*{Beyond our model}
	We show that the fittest individuals in the genealogy drive the adaptation in an asexual population facing changing environment, a feature already observed in other evolutionary context~\citep{HerRedWag02,neher2013genealogies,RouCof07}. Beyond this qualitative information, we quantify this phenomenon in specific contexts (diffusive approximation and small mutation regime). 
	In particular, in the regime of small mutations, we show that the structure of Hamilton Jacobi equations carries some information on the genealogies resulting from adaptation to this changing environment. 
	This opens a broad range of applications for the method presented in this work, since this Hamilton-Jacobi approach has been used extensively in different ecological models, 
	and, in particular, when mutations may have a large effect, for instance 
	with a kernel $K$ that does not satisfy the  assumption \cref{ass K} of being exponentially bounded. 
	Recently,  \cite{bghp18,mirrahimi2020singular} 
	found the (non-stationary) Hamilton-Jacobi equation when $\sigma=0$, in the case of 
	a broad range of fat tailed kernels. 
	However, in that case, the dual Lagrangian point of view of the equation is no longer valid, since the Hamiltonian $H$, in \cref{limit U}, is no longer well defined. However, our investigation of genealogies, based on neutral fractions, can still apply, and would provide a description of ancestral lineages in this context.
	From a mathematical point of view, we mainly use the linearity of $\mathcal B$ and spectral results about the linearized operator $\mathcal L$ around the equilibrium $F$. These properties hold true in many models, see \cite{mischlerscher}. 

	In our model, we assume that the environmental change equally affects all the individuals in the population. 
	However, habitats may differ between locations.  
	Recently, many works have focused on the interaction of two populations living in two different habitats, \textit{e.g.} \cite{hamel2020adaptation,mirrahimi2020evolution}. Each habitat favors a different optimal trait. Individuals can move between habitats and are subject to natural selection. In this scenario, migration tends to shift the trait distribution in each habitat towards the optimal trait of the other habitat. In particular, polymorphism can appear in a habitat provided migration is strong enough. 
	Investigating ancestral lineages in this context could yield valuable insights on this phenomenon. 
	Recently, \cite{GarLaf20} have extended the notion of neutral fractions to metapopulation model. Thus, using our notion of ancestral process based on neutral fractions can be an efficient tool to investigate the genealogy of a metapopulation located in different areas.

	An other important issue in spatial ecology is the evolution of populations undergoing range expansions. 
	Over the last decades, several theoretical works have focused on models of adaptation with  continuous space and trait variable, such as the ``cane toad equation'' to tackle these issues~\cite{benichou2012front,bouin2017super}. 
	In expanding population, the fittest individuals are at the front of the propagation range, a phenomenon called 'spatial sorting'. Up to our knowledge, the mechanisms underlying this phenomenon are poorly understood. Recently, \cite{canetoadslows} showed that a non local competition term can slow down the acceleration of the front, conversely to previous formal intuitive  results of~\cite{bouin2012invasion}. Thus, the investigation of the genealogy of individuals at the leading edge of the front, using our methodology, might provide new insights on this issue.
	
	In a related work, \cite{etheridge_genealogies_2020} studied a spatial Moran process modeling an expanding population with a strong Allee effect (corresponding to a bistable reaction-diffusion equation). 
	Using  neutral markers, they showed  that the genealogy of individuals sampled near the position of the traveling front is asymptotically close to the classical Kingman coalescent. 
	This result is in stark contrast with classical results on pulled waves (as in the stochastic Fisher-KPP equation) where the genealogy of individuals sampled at the leading edge of the front is believed to follow 
	the Bolthausen-Sznitzman coalescent, in which multiple lineages can merge at the same time \cite{berestycki_genealogy_2013-1,Brunet2007,desai2013genetic,neher2013genealogies}. 
	As it happens, in order to obtain their result, \cite{etheridge_genealogies_2020} use the fact that ancestral lineages of individuals sampled near the front approach a diffusion process which admits a stationary distribution in the moving frame centred on the front position, a fact that is reminiscent of what we obtain here.
	
	In the present paper, we focus on asexual population while many species reproduce sexually. 
	When this is the case, each individual has two parents and the reproduction can be described using the infinitesimal model~\cite{BarEthVeb17,Fisher1918,Tur17,TurBar94} which takes the following form in our framework:
	\begin{align}\label{sexop}
	\mathcal{B} (f)(x) :=  \dfrac{1}{\sigma \sqrt{\pi}}   \iint_{\mathbb{R}^{2}}  \exp\left[ -  \dfrac{1}{\sigma^2} \left( x - \dfrac{ x_1 +   x_2}2 \right)^2 \right ]   f( x_1)\dfrac{   f(  x_2)}{ \int_{\R}   f(  x_2')\, d  x_2'}\, d  x_1 d  x_2.
	\end{align}
	This operator states that individuals of trait $x_1$ and $x_2$ give birth to an individual whose trait $x$ is drawn from a Gaussian distribution centered at the mean of the traits of its parents, $(x_1+x_2)/2$, and with variance $\sigma^2.$ 
	Thus in this context, the genealogy of individuals becomes a binary tree which tracks the whole pedigree of each individual. 
	The analysis of ancestral lineages thus becomes more intricate, since the number of genealogical ancestors can quickly reach the size of the whole population.
	However,  recent asymptotic studies~\cite{spectralsex,timespectralsex} have laid the groundwork for the adaptation of the  neutral fractions framework to this setting.

	\subsection*{Stochastic framework, coalescence}$ $\\
	As already mentioned, the deterministic models~\cref{eq:f z} and~\cref{eq:TW_c_diff} can be obtained as large population limits  of stochastic individual-based models described in \Cref{app:IBM}~\cite{champagnat2006unifying,champagnat2007individual}. 
	In this stochastic context, the neutral fractions approach has also been used to understand how lineages coalesce back in time. For instance, recently \cite{billiard_stochastic_2015} studied neutral markers in the background of a trait-substitution sequence in the adaptive dynamics regime, see also \cite{etheridge_genealogies_2020} for genealogies in bistable traveling wave.
	
	Although our deterministic model cannot track the microscopic aspect of the lineages because the coalescence events become increasingly rare as the size of the population tends to infinity, we can still provide some heuristics on the timescale of coalescence in the genealogy of individuals.
	For instance, \Cref{fig lineage asex} suggests that lineage are less likely to coalesce before they reach a trait close to the optimal trait. 
	Moreover, in our setting, the growth rate of individuals is maximal near the moving optimum, and becomes negative far ahead of the optimum, so that lineages do not escape in the tip of the front (as they do in pulled fronts).
	This means that the travelling pulse of \cref{eq:TW_c_diff} is ``pushed'' and should behave similarly to bistable waves, as studied in \cite{etheridge_genealogies_2020}.
	As a result, most coalescence events should take place near the optimum, and no single individual is likely to quickly produce a large progeny (of a size comparable to the total population).
	We can thus conjecture that, as in \cite{etheridge_genealogies_2020}, the genealogy of a finie sample of individuals follows Kingman's coalescent. 

	
	In addition, using some of our results, we can provide heuristics about the mean coalescence time of two lineages $T_2$. First, our analytic approximations \cref{mean var quadex} and \cref{gamma quad approx} on the dynamics of the mean of the ancestral lineage provides a good approximation on the characteristic time $T_0$ before which two lineages reach the optimal value $0$. From our diffusive approximation we get 
	\[
	T_{0}\approx \dfrac1{\sigma\sqrt{\beta}}
	\]
	In the case of more general operators, based on the approximation formula~\cref{gamma approx} in the small variance regime, we expect
	\[
	T_{0}\approx \dfrac 1{\Var(F)},
	\]
	where $\Var(F)$ is the variance of the phenotypic distribution at equilibrium.
	After this time delay $T_0$, the coalescent time $T_2$ seems to follow an exponential distribution (see \Cref{fig_app:time_coalescent}). However, the parameter of the exponential time depends on the speed of the changing environment. If $c$ is small ($0<c\leq \sigma$), our simulations suggest an exponential rate $T_0$ while if $c$ is large ($c>\sigma$) the exponential rate seems to depend on the distribution of the ancestral process and $N$ size scale of the population (see \Cref{fig_app:time_coalescent}). When the changing speed $c$ is small, the amount of individuals at the optimal trait $ct$ is large compared to the size of the population. Thus, using the heuristic of~\cite{etheridge_genealogies_2020,neher2013genealogies} with our \cref{long time asymp}, we expect that the exponential rate does not truly depends on the population size $N$ and it should be $T_0.$ However, when the speed increases, the population size at the optimal trait is low compared to the size of the population, thus the exponential rate should depend on it. Using the heuristic of~\cite{etheridge_genealogies_2020}, we suggest that the exponential rate should be 
	$$ \dfrac{1}{2\sqrt{N}} \int_\R \dfrac{1}{F(z)} \dfrac{(F(z)\varphi(z))^2}{\left(\int_\R F(z) \varphi(z) dz \right)^2} dz.$$
	These arguments are all in the preliminary stage of a larger work, and need to be further investigated. 
	
	

	
	%
	%
\begin{small}	\bibliographystyle{apalike}
	\bibliography{Biblio,bib_raph}
	\end{small}

	\begin{appendices}
		\crefalias{subsection}{appendix}
		
		
		\counterwithin{figure}{section} 
		\section{The regime of small variance : further developments}
		\subsection{Derivation of the Hamilton-Jacobi equation}\label{sec HJ derivation}$ $\\
		In this section we explain formally how to derive a Hamilton-Jacobi equation from the pulse characterization  \cref{eq:equilibrium c}, in the regime of small mutations. 
		From the work of \cite{main}, we know that the regime of small mutation effects can be described by the following scaling factor:
		\begin{align}\label{def eps}
		\e := \sqrt{\sigma^2 \frac{\alpha}{\beta}}, \text{with } \alpha :=  \mu '' (0)>0.
		\end{align} 
		When $\e \to 0$ our equation~\cref{eq:equilibrium c} converges to a Hamilton-Jacobi equation which allows us to use the rigorous results of \cite{lorz2011dirac}. This limiting regime captures the \emph{weak selection regime}, in the sense that either the variance vanishes, or the  selection is weak compared to birth: 
		$$ \sigma \ll 1 \quad \text{ or } \quad  \frac \alpha \beta  \ll 1. $$
		We show here how the pulse $F$ concentrates around a mean trait value, $z^*$,  when $\sigma \to 0$. 
		Recall that $F$ is a solution of the integro-differential equation \cref{eq:equilibrium c}, which we rewrite here: 
		\begin{align}\label{spec asex}
		\lambda_\sigma F_\sigma(z) - c  F'_\sigma(z)  + \mu(z) F_\sigma(z) = \dfrac{\beta}{\sigma} \int_\R  K\left(\dfrac{z-z'}\sigma\right)F_\sigma(z')\, dz' .
		\end{align}
		%
		When the variance of the mutation kernel vanishes, that is $\sigma \to 0$,  we may expect the function $F$ to concentrate around a specific trait, guided by selection, see \cite{barles2009,lorz2011dirac} for instance. To capture this phenomenon, we use the following Hopf-Cole transform, identical to \cref{limit U}:
		\begin{align}\label{hopf asex}
		F_\sigma(z) = \exp \left( - \frac{U_\sigma(z)}{\sigma}\right).
		\end{align}
		The scaled quantity $U_\sigma$ then satisfies the following equation::
		\begin{align}\label{eq:spec asex}
		\lambda_\sigma  + \frac{c}{\sigma}  U'_\sigma(z)  + \mu(z) = \dfrac{\beta}{\sigma} \int_\R  K\left(\dfrac{z-z'}\sigma\right) \exp \left( \frac{U_\sigma(z)-U_\sigma(z')}{\sigma}\right) \, dz' .
		\end{align}
		To avoid degeneration of terms in the equations when $\sigma\to0$, 
		we rescale the speed of environmental change by $\sigma$,  as in \cref{scaled c}:
		\begin{align}\label{scaling c asex}
		c := \sigma c'.
		\end{align}
		In order to find the limit equation, we use the following Taylor expansion of the exponential term inside the integral:
		\begin{align*}
		\forall z'\in \R, \quad  \exp \left( \frac{U_\sigma(z)-U_\sigma(z')}{\sigma}\right) \approx \exp \left(-\frac{(z-z')}{\sigma} \, U'_\sigma(z)\right).
		\end{align*}
		Plugging this approximation into  the integral term of \cref{eq:spec asex}, with an affine change of variable,  we obtain  formally when $\sigma \to 0$:
		\begin{align*}
		\lambda_0  + c' U'_0(z)  + \mu(z) = \beta  \int_\R  K(y) \exp \Big( y\,  U'_0(z) \Big) \, dy .
		\end{align*}
		In the following, we omit the index $0$, as in \cref{limit U}. 
		We then recover equation \cref{spec Ham}. Let us observe that the limit equation is well defined thanks to assumption~\cref{ass K} on the exponential decay of $K$. It also guarantees the finitess of  the Hamiltonian $H$  in \cref{limit U} for all $p \in  \R$.
		

		\subsection{Lagrangian point of view : qualitative formulas}\label{sec Lag pov}$ $\\
		The Hamilton-Jacobi equation \cref{spec Ham} provides analytical formula as well as qualitative behavior on $U$ and $\lambda$, which are expected to be a good approximation when $\sigma$ is small. Let $\mu_0:=\mu(0)$ be the minimum mortality rate.  
		
		We first start with the following formula on  $\lambda$:
		\begin{align}\label{eqint asympGamma4}
		\lambda= \beta -\mu_0- \beta  L \left( \frac{c'}{ \beta}\right).
		\end{align}
		where $L$ is the Lagrangian associated to the Hamiltonian $H$ and related to the mutation kernel $K$, see \cref{def Lag}. 
		The quantity $\beta-\mu_0$ corresponds to the intrinsic growth rate (fitness) of the population while the additional quantity $\beta L( c'/  \beta ) $ measures the lag load induced by the changing environment.
		
		A short argument for \cref{eqint asympGamma4} consists in assuming that the asymptotic behavior of $\Gamma$ stated in \cref{lim gamma} holds true. Then, plugging this into \cref{U0 Gamma}, it  prescribes the value of $\lambda$ such that $U$ does not take infinite values: 
		\begin{align*}
		\beta L\left(\frac{c'}{\beta}\right) - \beta + \mu(0) + \lambda = 0.
		\end{align*}
		Since $\mu(0)=\mu_0$, this formula coincides with \cref{eqint asympGamma4}. 
		
		The proof of this analytical formula relies on  convex analysis methods (see \cite{main}).  First, the function $p\mapsto c'p-\beta H(p)$ admits a maximum value : $\beta L(c'/\beta)$. Indeed, the functions $H$ and $L$ have reciprocal derivatives functions: $\p_p H \circ \p_v L = Id$.  Adding this maximum value on each side of the Hamilton-Jacobi equation \cref{spec Ham}, we obtain
		\begin{align}\label{eqint asympGamma3}
		\mu(z)-\mu_0 + \left[  \lambda-\beta + \mu_0  + \beta L \left( \frac{c'}{ \beta}\right) \right]= \beta H( U'(z)) -c'  U'(z) +  \beta L \left( \frac{c'}{ \beta}\right).
		\end{align}
		We claim that  on the  left hand side of this relationship, the term between brackets  must vanish. Otherwise, it would lead to a contradiction as follows. 
		
		Indeed, the function $\mathcal C : p\mapsto \beta H(p) -  c'p + \beta L(c'/\beta)$ is convex, nonnegative and reaches zero from the properties  of the Hamiltonian $H$ and the Lagrangian $L$. If the term between brackets does not vanish in \cref{eqint asympGamma3}, the function $z \mapsto  \mathcal C \circ U' $ takes only (strictly)  positive values. As a consequence, $ U'$ only takes values in one of the two branches of the convex function $\mathcal C$, and on each of these branches, $\mathcal C$  is invertible. Therefore, for each $z \in \bar \R$, we can invert the relationship \cref{eqint asympGamma3} to deduce the value taken by $U'(z)$. From \cref{hyp growth m}, $\mu$ has the same infinite value at $\pm \infty$. Inverting \cref{eqint asympGamma3} for $z = \pm \infty$ yields $$\lim_{z \to  + \infty} U'(z)= \lim_{z \to  - \infty} U'(z).$$
		This is in contradiction with the  assumption $U(\pm \infty)= + \infty$, or equivalently, $F(\pm \infty)=0$, \textit{i.e.} the population density vanishes at infinity. Therefore, the term between brackets in \cref{eqint asympGamma3} vanishes, exactly as in the desired formula \cref{eqint asympGamma4}. 
		
		As a side,   $z=z^\ast$ leads to a formula that dictates the  position of the dominant trait  $z^\ast$:
		\begin{align*}
		\lambda + \mu (z^\ast) =\beta.
		\end{align*}
		Combined with \cref{eqint asympGamma4}, we find
		\begin{align}\label{loc zast asex}
		\mu(z^\ast)  =\mu_0 + \beta  L \left( \frac{c'}{ \beta}\right).
		\end{align}
As a matter of fact those formulas are consistent with those in \cite{main},
		\begin{align*}
		\lambda \approx \beta- \mu_0- \beta L\left(\frac{c'}{ \beta} \right) + O(\e).\end{align*}
		They further show how to obtain more accurate expansions,  up to an arbitrary  order (in $\e$ defined in \cref{def eps}), and compute explicitly the following corrector terms. 
		
		In addition,  we know  from \cref{eq:rhob} that $\lambda$ is a measure of the size of the population. Thus, this size should remain positive which gives us a critical threshold for the speed $c'$ so that the population does not go extinct:
		\begin{align*}
		c_{max} =  \sigma \beta \, L^{-1} \left( \frac {\beta-\mu_0} \beta \right).
		\end{align*}
		Moreover, we can check by integrating~\cref{eq:equilibrium c} that 
		\[
		\lambda = \int_\R\dfrac{\mu(x)F(x)}{\int_\R F(x')dx'}dx
		\]
		which corresponds to the mean fitness of the population, or its mean intrinsic rate of increase.

		\subsection{Long time behavior of $\Gamma_z$}\label{sec lim Gamma}
		In this part, we  show a somehow stronger version of  \cref{long time asymp} : the typical ancestral lineages $\Gamma_z$ converges towards $0$, the optimal trait, when $s$ goes to infinity. In the regime of small variance, this result concerns  the ODE \cref{edo Gamma}. 
		
		Let $z_0$ be a steady state  of \cref{edo Gamma}. Then
		\begin{equation}\label{eqint Asymp7}
		0 = -c'+\beta \, \partial_p H(U'(z_0)).
		\end{equation}
		Let  $p_0 := U'(z_0)$. Then, from \cref{eqint Asymp7}, $p_0$ is a critical  point of the function $p\mapsto  - c' \, p +   \beta \, H(p)$.  Since this is a convex function, we deduce that
		\begin{align*}
		p_0 \in \argmin_{p} \Big(-c' p +  \beta \,  H(p)\Big).
		\end{align*}
		Therefore, by definition of the Lagrangian function $L$ in \cref{def Lag}
		\begin{align}\label{eqint Asymp8}
		- \beta L\left(\frac{c'}\beta\right)  = -c'p_0 +\beta \, H(p_0).
		\end{align}
		On the other hand, by evaluating the Hamilton Jacobi equation \cref{spec Ham} in $z = z_0$, one finds
		\begin{align*}
		\mu(z_0) + \lambda   =  \beta +\beta \, H(p_0)  -c'p_0.
		\end{align*}
		Plugging in \cref{eqint Asymp8}
		\begin{align}
		\mu(z_0)   = \beta- \lambda +\beta L\left(\frac{c'}\beta\right). 
		\end{align}
		Finally, thanks to the formula \cref{eqint asympGamma4} for $\lambda$, we get
		\begin{align}
		\mu(z_0)   = \mu_0.
		\end{align}
		According to our assumptions on $\mu$ stated in \cref{hyp growth m}, $0$ is the only global point of minimum of the convex function $\mu$, therefore $z_0=0$. We conclude that $0$ is the unique steady  state of the ODE \cref{edo Gamma}. Moreover, it is established in \cite{lorz2011dirac} that under the assumptions \cref{ass K}, $U$ is a convex function. Therefore, the flow in \cref{edo Gamma} is an increasing function, and it is straightforward that, $\Gamma_z$ converges to this unique steady state of \cref{edo Gamma}: $$\Gamma_z(s) \xrightarrow[s \to \infty]{} 0.$$

		\subsection{Approximation of the mean of the ancestral lineages}\label{sec approx app}$ $\\
		Working further on the Hamilton Jacobi  equation~\cref{spec Ham}, we can make an analytic approximation of $\Gamma_z$, the typical lineage. Let us first 
		differentiate \cref{spec Ham} with respect to $z$, and then divide by $U''$ on each side. We obtain for  all $z\in\R$,
		\begin{align}\label{eq approx}
		-c'+ \beta \, \p_p H \big( U' (z)\big) = \dfrac{\mu '(z)}{U'' (z)}.
		\end{align}
		Now, our idea is to link $U''(z)$ with the variance of $F$.  Using results from  \cite{main}, we obtain the following approximation of the variance of $F$ at the leading order in $\sigma$:
		\begin{align}\label{var sigma0}
		\Var(F) = \frac \sigma {U''(z^*)} + o(\sigma),
		\end{align}
		where  $z^*$ is the dominant trait in our population. This approximation comes from a Taylor expansion (with Laplace's method) of the integrals defining the variance:
		\begin{align*}
		\Var(F) = \frac{\ds \int_{\mathbb{R}} (z-z^\ast)^2 \exp \left( - \frac {U(z)}{\sigma} \right) dz}{\ds \int_{\mathbb{R}} \exp \left( - \frac {U(z)} \sigma \right)dz }.
		\end{align*}
		In addition, we make the following rough approximation, valid if $z$ is close to $z^*$: $U''(z)\approx U''(z^*)$. Plugging these approximations into \cref{eq approx}, we find that 
		\[
		-c'+ \beta \, \p_p H \big(  U'(z)\big) \approx \dfrac{\mu'(z)\Var(F)}{\sigma} \ \hbox{ for } \ z \ \hbox{ close to } \ z^*.
		\]
		The ODE \cref{edo Gamma} satisfied by $\Gamma_z$ then becomes:
		\begin{align}\label{edo Gamma_approx}
		\dot \Gamma_z(s)& = - \dfrac{\Var(F)}{\sigma} \mu '\big(\Gamma_z(s)\big) ,\\
		\nonumber \Gamma_z(0) & =z .
		\end{align}
		In particular, \textbf{if the selection function is quadratic},  
		$\mu (z)= \mu_0 + z^2/2$, the solution of \cref{edo Gamma_approx} is  precisely
		\begin{equation}\label{gamma approx}
		\Gamma_z(s) = z \exp\left(-\dfrac{\Var(F)}{\sigma}\,s\right).
		\end{equation}
		Moreover, on this case, we can derive an explicit formula for the variance from the HJ equation.
		More precisely, we compute $U''(z^\ast)$  from \cref{spec Ham} which provides the following formula from \cref{var sigma0}:
		\begin{align*}
		\Var (F)= -\frac{c' \sigma }{\mu '(z^\ast)}.
		\end{align*}
		Since $\mu$ is quadratic and  $z^\ast$ satisfies \cref{loc zast asex}, we obtain
		\begin{align}\label{var quad}
		\Var(F) = \dfrac{c'\sigma}{\sqrt{2\beta}\sqrt{L\left(\dfrac{c'}{\beta}\right)}}. 
		\end{align}
		Finally, we get the following approximation for the mean trajectories of the ancestral process: 
		\begin{equation}\label{gamma quad approx}
		\Gamma_z(s) = z \exp\left(-\frac{c'\,s}{\sqrt{2\beta}\sqrt{L\left(c'/\beta\right)}}  \right).
		\end{equation}
		This formula is explicit, since it only depends on the mutation kernel through the Lagrangian $L$. However this approximation only applies to the case of  a quadratic selection and for traits close to the dominant trait. But we can check from numerical simulations, that  this approximation is quite robust (see \Cref{fig:mean_approx}).
		
		
		\section{Numerical methods}\label{sec appendix num}
		\subsection{Ancestral lineages in the stochastic model}\label{app:IBM}$ $\\
		We detail in this section the  numerical simulations how we deal with the simulations of the individual based model and the lineages. 
		The population at each time $t$ is made of a number $N(t)$ of  alive individuals.   Let us consider an individual, denoted $i$, that is alive at time $t$, with $1\leq i \leq N(t)$. It has the trait $z_i\in \R$. The first event  for this individual is one of the following :  
		\begin{itemize}[label=$\triangleright$]
			\item \textbf{Birth of a descendant:} it happens at an uniform rate among individuals $r_B^i= \beta$. 
			\item  \textbf{Death of the individual $i$:} it  is decomposed in two separate events: \begin{enumerate}
				\item[$(a)$] \textbf{ Death by selection} The individual may die because its phenotype is ill-adapted in the phenotype landscape. This happens at the rate:
				\begin{equation*}
				r_{Ds}^i = \mu(z_i-ct).
				\end{equation*}
				\item[$(b)$]\textbf{ Death by competition} Alternatively, an individual may die because of the density dependence in the population, at a rate that that depends on the total size size of population at  time $t$,  and on the carrying capacity $N$
				\begin{equation*}
				r_{Ddd} = \ds \sum_{j=1}^{N(t)}\frac{1}{N}=\frac{N(t)}{N}.
				\end{equation*} 
			\end{enumerate}
		\end{itemize}
		\textbf{Next event : incrementation of the time step} The time step is the smallest time for all individuals to go through one of the previous steps. Thanks to the  Markov property, each event occurs following an exponential law of parameter:  
		\begin{align}\label{next event}
		dt \sim  \min_{1 \leq i \leq N(t) }\mathcal{E}(r_B^i+r_{Ds}^i + r_{Ddd}).
		\end{align}
		By the memory loss property of the exponential distribution, $dt$ also  can be drawn from an exponential law which rate is the sum of the rates of all the independent events: $$dt \sim \mathcal{E}\left(\ds \sum_i (r_B^i+r_{Ds}^i + r_{Ddd})\right).$$
		%
		
		\textbf{Update of the population:} Once the next event is decided, according to the law \cref{next event}, the population at time $t+dt$ is deduced by either adding the individual that was born ($N(t+dt)=N(t)+1$) or subtracting the one that died ($N(t+dt)=N(t)-1$). In the case of a birth event, the trait of the offspring is drawn according to the operator $\mathcal B$ in \cref{eq:asexual B}:
		$$ z_\text{offspring} = z_\text{parent} + \sigma dK.$$
		We repeat all the steps until reaching the  final time of simulation. 
		Numerically, this model has a very high computational cost, because it needs a relatively high number of individuals to approach the deterministic model given by \cref{eq:f z}. 
		As a consequence, we performed the simulations using an approximating model, by first fixing $dt$ to a small but deterministic value. Then, for each individual, we draw a time of birth following the law $\mathcal{E}(\beta)$ and a time of death following the law $\mathcal{E}(\mu(z_i)+N(t)/N)$. Then  we simply count which individuals led to a reproduction event and which died on the time-window $\left[t,t+dt\right]$. This amounts to the supposition that on this time interval, individuals cannot reproduce more than once.

\begin{figure}[H]
			\begin{center}
				\includegraphics[scale=0.35]{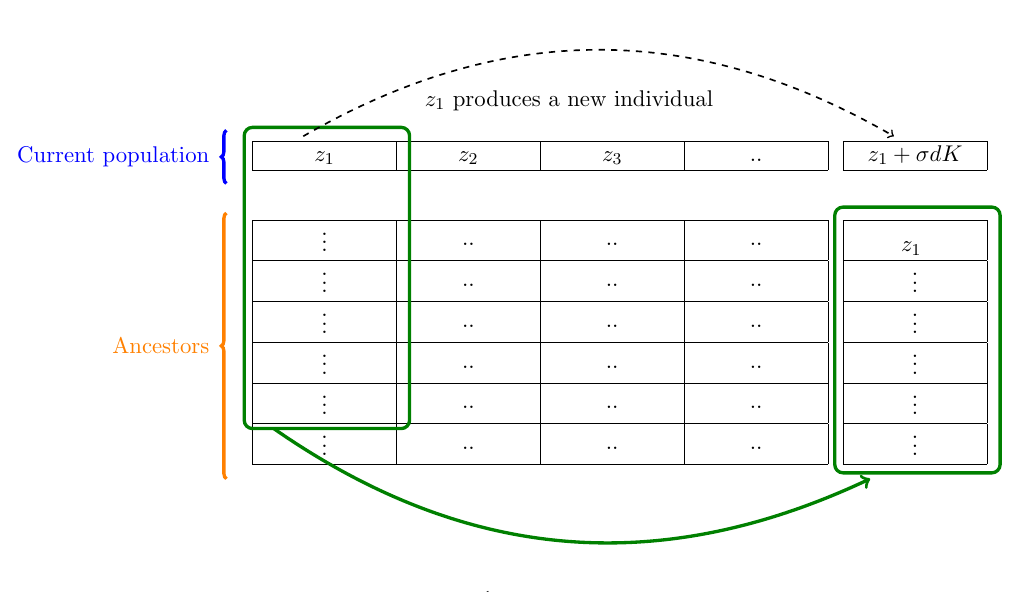}
			\end{center}
			\caption{\small Keeping track of the lineages : a schematic representation of the ancestral lineages algorithm described in~\Cref{app:IBM}}\label{fig tabnum}
		\end{figure}

		Finally, to follow the lineage of individuals, we create a huge matrix at the start of the simulation. We will stock the lineage of every individual in this matrix, filling it progressively.  Every time an individual appears, its lineage is similar to the one of its parent, translated by one generation.  The numerical procedure works as described in \Cref{fig tabnum}, where different  columns correspond to different individuals, and each line corresponds to a new generation.
		
\begin{figure}
			\begin{center}
				\includegraphics[scale=0.45]{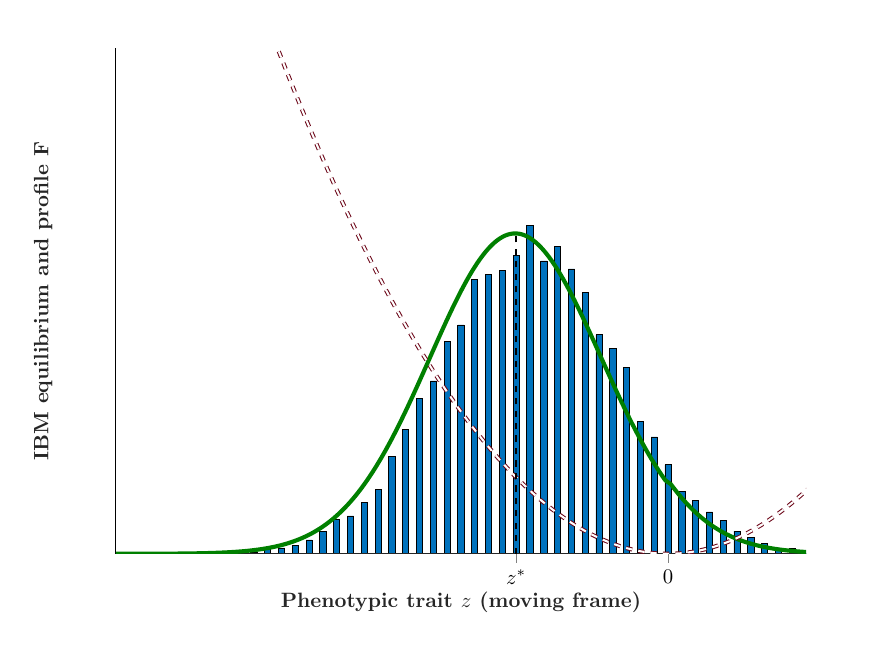}
				\caption{\small Population density profile $F$ in the moving frame obtained from deterministic model~\cref{eq:f z} (plain green curve) and from the IBM model with a number of individuals of order $2\cdot10^4$ (blue histogram). The selection function $\mu$ is quadratic with minimum at $z=0$ (red dashed curve). The value $z^*$ corresponds to the mean phenotypic trait of the deterministic density.}\label{fig profiles asex}
			\end{center}
		\end{figure}

		This algorithm led to \Cref{fig profiles asex,fig lineage asex mean var bis} with the following parameters: 
		\begin{align*}
		\alpha = 2, \quad \beta = 2, \quad \sigma = 0.1, \quad N=20000, \quad c=0.2\, .
		\end{align*}

	\newpage
	\subsection{Mean traits along the lineages}
		$ $
		We extend the simulations of \Cref{fig lineage asex mean var bis} to different scenarios corresponding to various  mutational variances $\sigma$ and speeds $c$. More precisely, for two different couples of mutational variance and speed, we compare the mean and the variance of the stochastic lineages with the mean and the variance of the ancestral process defined by our PDE model. We show that our deterministic model provides agreed with the individual based model. In addition, we see that the variance of the ancestral process increases with the mutational variance $\sigma$ while it slightly increases with the speed $c$.  
			\begin{figure}[H]
			\subfigure[$\sigma=0.05$ and $c=0.025$]{\includegraphics[scale=0.5]{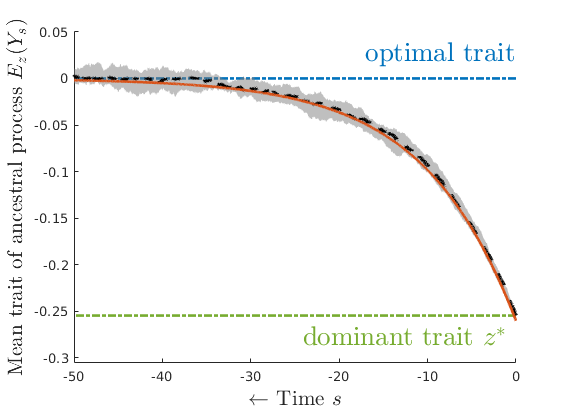} \includegraphics[scale=0.5]{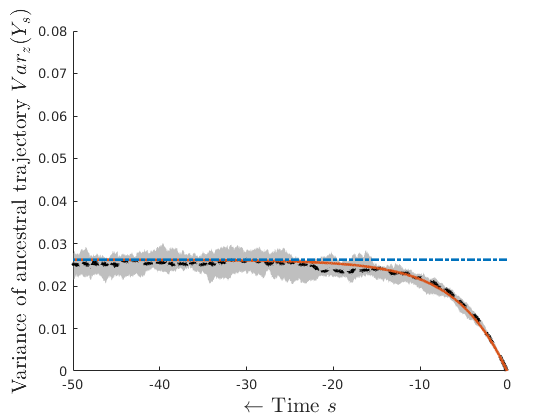}}
			\subfigure[$\sigma=0.1$ and $c=0.05$]{\includegraphics[scale=0.5]{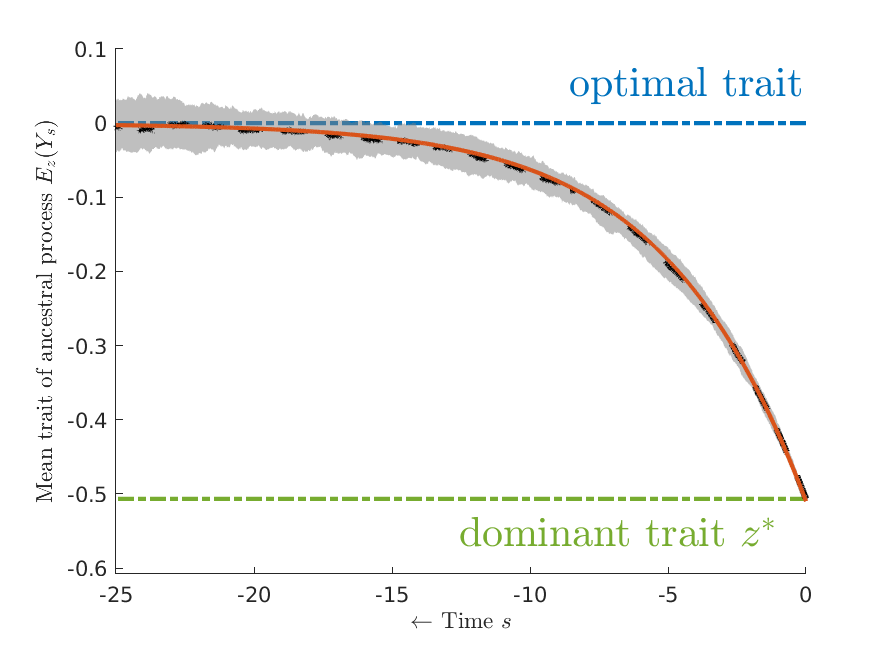} \includegraphics[scale=0.5]{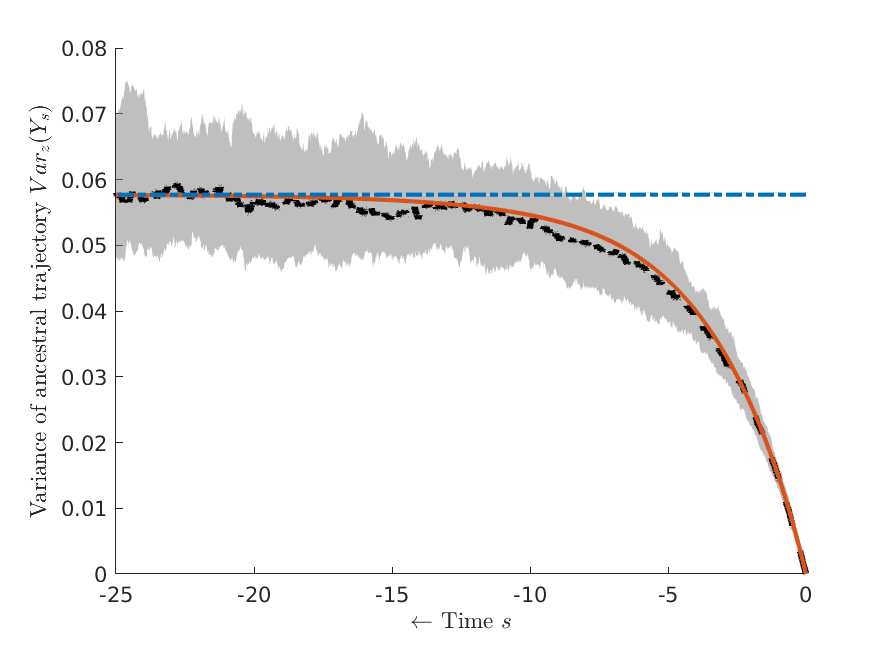}}
			\caption{ \small \textbf{Mean (left) and variance (right) of the ancestral process $Y_s$ for different set of parameters}: red curves corresponds to the model~\cref{dual process} and the black dashed curves correspond to the IBM model averaged over $50$ replicates and the gray regions corresponds to the $5\%$ and $95\%$ confident interval.   
				On the left panel, the horizontal blue lines are the optimal trait in the moving frame, $0$. The green lines are the dominant trait of the equilibrium $F$, denoted $z^\ast$. On the right panel, the horizontal cyan line corresponds to the asymptotic variance given by the deterministic model see~\cref{long time asymp}.
			}\label{fig app lineage asex mean var}
		\end{figure}

		\newpage
		\subsection{Coalescent time}
		
\begin{figure}[h]
			\begin{center}
				\subfigure[$\sigma=0.1$ and $c=0.1$]{\includegraphics[scale=0.5]{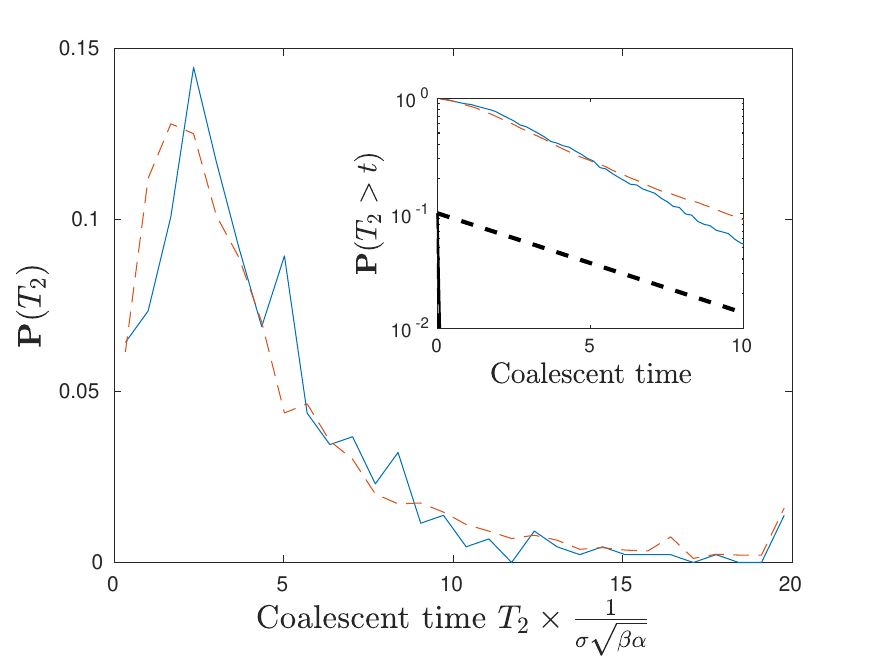}}
				\subfigure[$\sigma=0.1$ and $c=0.2$]{\includegraphics[scale=0.5]{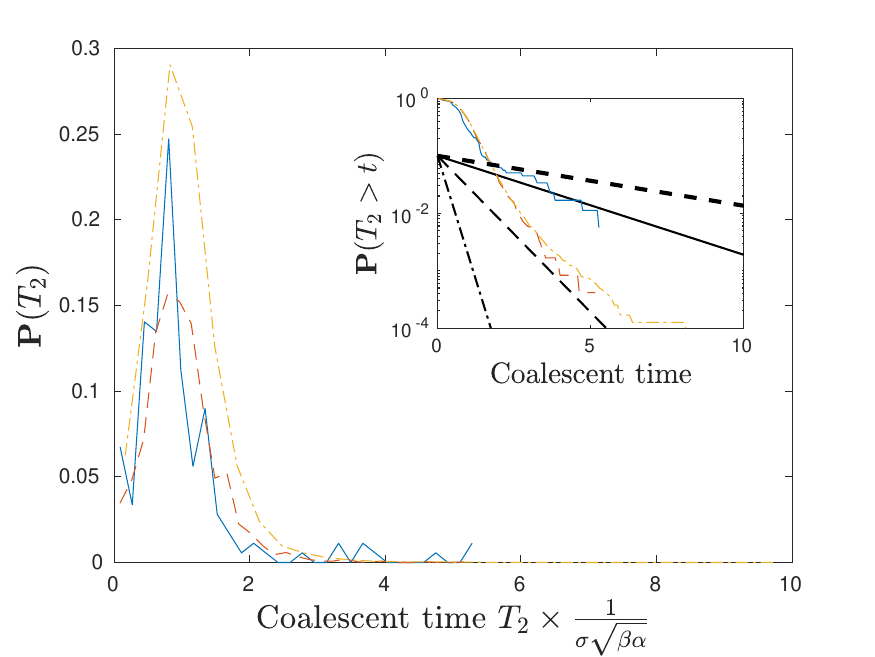}}
				\caption{\small Distribution of the time $T_2$ to the most recent common ancestor of two individuals for different speed $c$ (a) $c=0.1$ and b) $c=0.2$) and typical size $N$ of the population. The inset shows the complementary of the cumulative distribution function $\mathbf{P}(T_2>t)$ in the semilog coordinates.
				An exponential $\exp(-t/T_0)$ is indicated as a black dashed line. Different line style and colour corresponds to $N=10^4$ (plain blue), $N=10^5$ (dashed red) and $N=10^6$ (dot-dashed orange). The mutational variance is fixed to $\sigma^2=0.1$ and $\alpha=\beta=2$. 
					}\label{fig_app:time_coalescent}
			\end{center}
		\end{figure}

		We first investigate the time $T_2$ before two individuals lineages meet in the past. More precisely, for each individual in the population, we compute the minimal time such that its lineage coalesces with an other lineage in the population. From \Cref{fig_app:time_coalescent}, we can observe that after a time delay $T_0=1/(\sigma\sqrt{\alpha\beta})$, the time $T_2$ is exponentially distributed. The exponential distribution of $T_2$ is apparent from the inset in \Cref{fig_app:time_coalescent} which shows the complementary of the cumulative distribution function $\mathbf{P}(T_2>t)$.
		
		
\begin{figure}[h]
			\begin{center}
				\subfigure[$\sigma=0.1$ and $c=0.1$]{\includegraphics[scale=0.5]{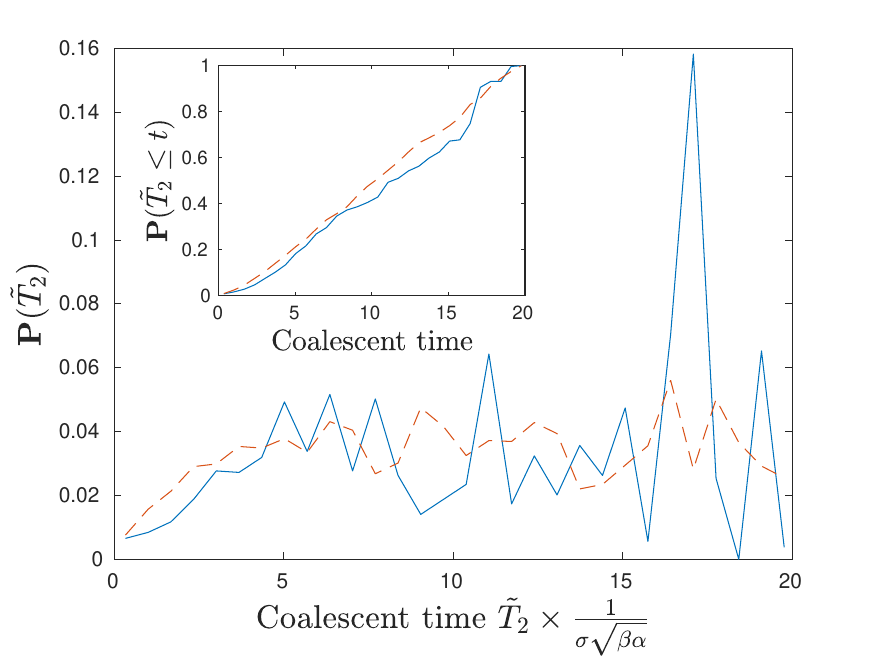}}
				\subfigure[$\sigma=0.1$ and $c=0.2$]{\includegraphics[scale=0.5]{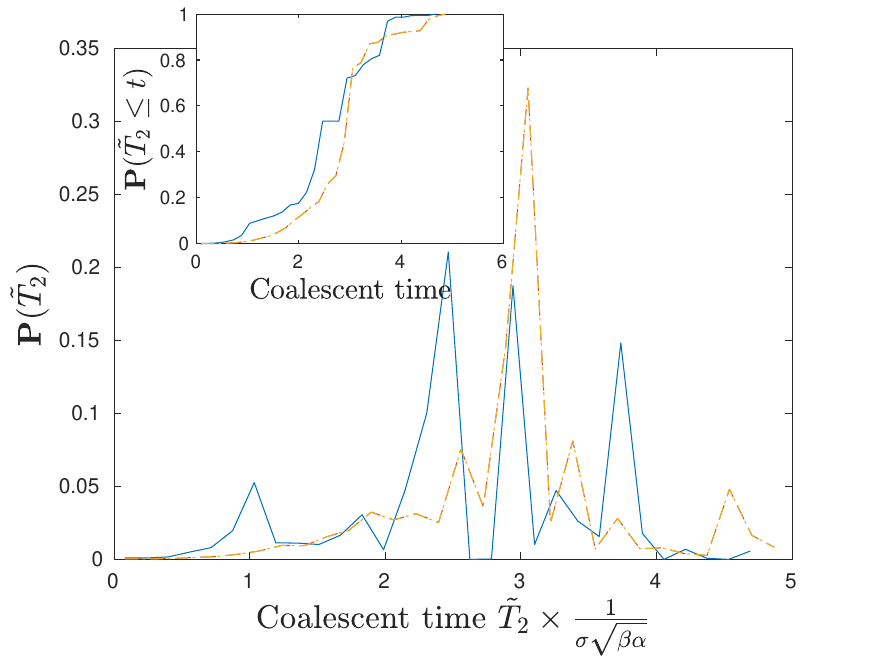}}
				\caption{\small Distribution of the time $\tilde T_2$ to the most recent common ancestor of two individuals sampled randomly in the population for different speed $c$ and typical size $N$ of the population. The inset shows the cumulative distribution function $\mathbf{P}(T_2\leq t)$. Different line style and colour corresponds to $N=10^4$ (plain blue), $N=10^5$ (dashed red) and $N=10^6$ (dot-dashed orange). The mutational variance is fixed to $\sigma^2=0.1$ and $\alpha=\beta=2$.}\label{fig_app:time_coalescent_tilde}
			\end{center}
		\end{figure}

		The time $T_2$ is different from the time $\tilde T_2$ before the lineages of two individuals sampled randomly in the population coalesce (see \Cref{fig_app:time_coalescent_tilde}). 

		\newpage
		\subsection{The PDE approach of \Cref{sec heur pde}}\label{app:IBM_vs_PDE}
		$ $
		
\begin{figure}[H]
			\begin{center}
				\subfigure[$\sigma=0.1$, $c=0.2$ and $N=10^4$]
				{\includegraphics[width=0.45\textwidth]{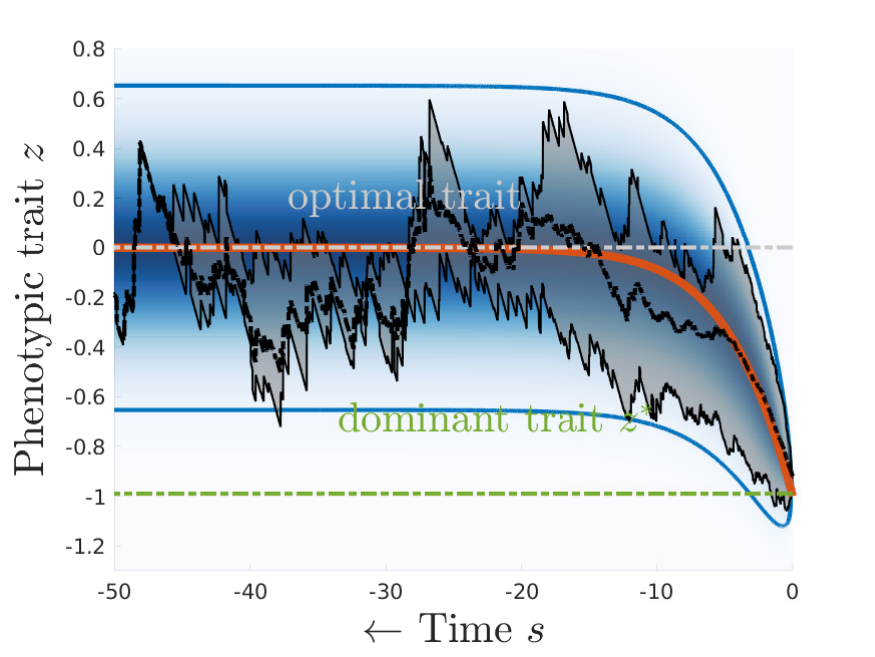}}
				\subfigure[$\sigma=0.1$, $c=0.1$ and $N=10^4$]
				{\includegraphics[width=0.45\textwidth]{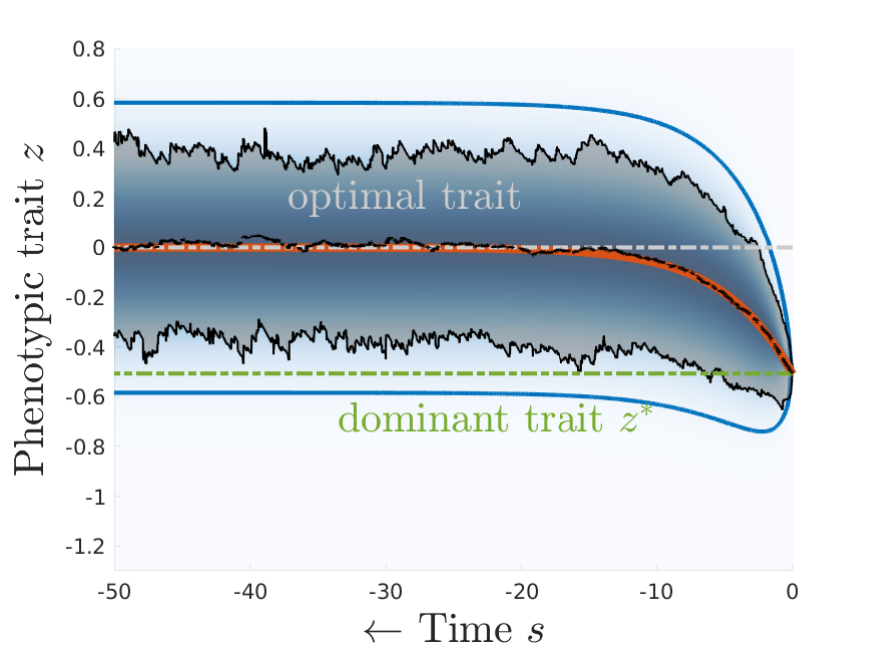}}
				\subfigure[$\sigma=0.1$, $c=0.2$ and $N=10^5$]
				{\includegraphics[width=0.45\textwidth]{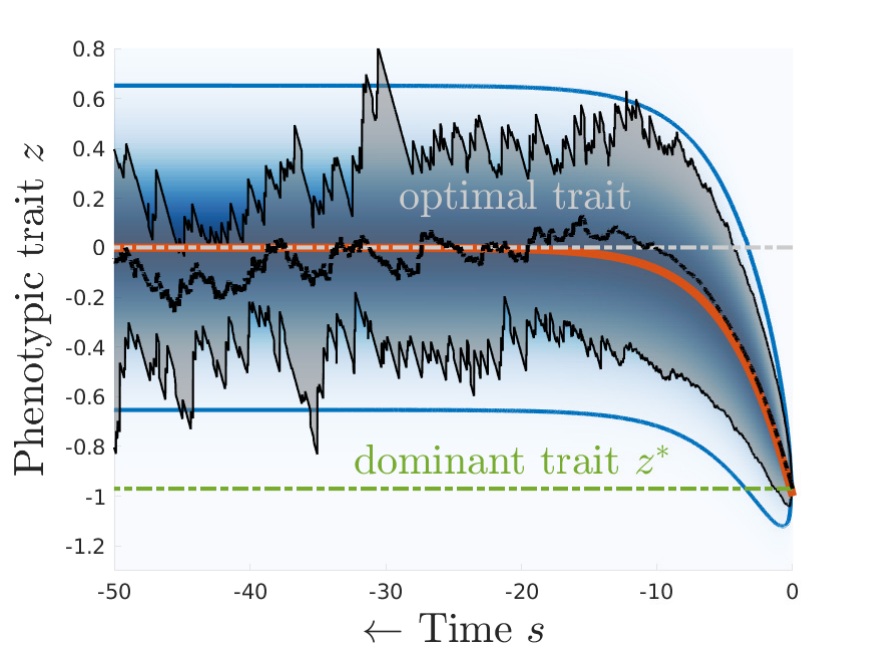}}
				\subfigure[$\sigma=0.1$, $c=0.1$ and $N=10^5$]
				{\includegraphics[width=0.45\textwidth]{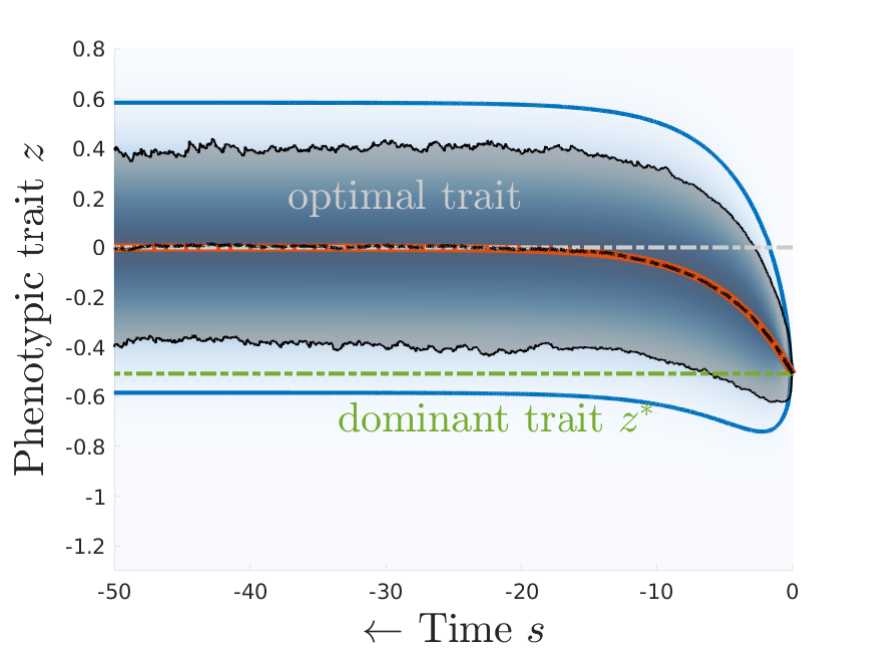}}
				\subfigure[$\sigma=0.1$, $c=0.2$ and $N=10^6$]
				{\includegraphics[width=0.45\textwidth]{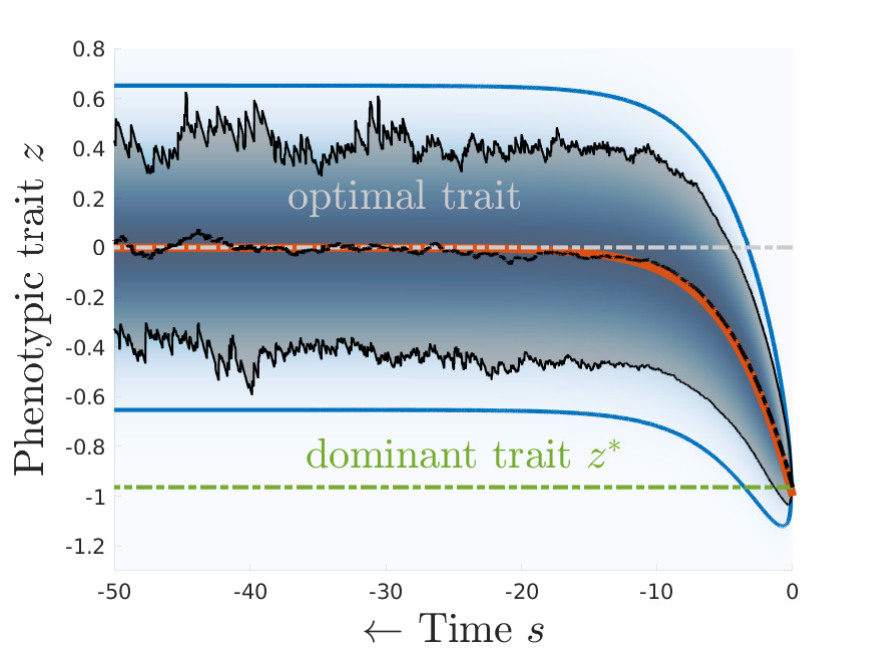}}
				\caption{Ancestral lineages over time of individuals with the dominant trait $z^*$: distribution of the ancestral process $Y_s$ (blue region) and the mean of the ancestral lineages of one realization of the IBM model (dashed black curve). The plain black curve and blue curve correspond to the confidence interval at $5\%$ around the mean for the IBM lineages (black) and the ancestral process (blue). The mean of the ancestral lineages corresponds to the dashed black curve and the mean of the ancestral process $E_z(Y_s)$ is represented by the plain red curve. The grey and green dot-dashed curves represent respectively the optimal and the dominant trait.  
				}\label{fig lineage asex distrib}
			\end{center}
		\end{figure}

		In this section we compare  the realizations of the IBM model with the ancestral process $Y_s$, obtained by the simulation of the fundamental solution defined by \cref{pde w}.  As mentioned in  \Cref{sec heur pde}, the distribution $w^z(s,y)$ of the ancestral process $Y_s$ starting from a trait $z$ is given by 
		\[
		w^z(s,y) = \dfrac{\upsilon^y(-s,y)}{F(z)}
		\]
		where $\upsilon^z$ satisfies
		\begin{equation}\label{eq_app:f_sub_dirac}
		\left\{ \begin{array}{ll}
		\p_t  \upsilon^y(t,z) & = \mathcal L(\upsilon^y(t,\cdot))(z),\\
		\upsilon^y(0,z) & = \delta(z-y)F(z).
		\end{array}\right.
		\end{equation}
		In order to solve numerically this equation, we replace the Dirac $\delta(z-y)$ with the characteristic function of a small interval around $y$ of the form $[y-dy,y+dy]$ as pictured in \Cref{fig neutral frac}. More precisely, we solve equation~\cref{eq_app:f_sub_dirac} on a finite interval of the form $[z_{min},z_{max}]$ with $z_{min} = -z_{max} = - 10$, and we add Dirichlet boundary condition. For the initial conditions, we evenly decompose the interval $[z_{min},z_{max}]$ into $n=200$ intervals $I_k$ of size $dy = (z_{max}-z_{min})/n$ and centered around $y_k$. The numerical solution of $\upsilon^k$ satisfying~\cref{eq_app:f_sub_dirac} starting from $\upsilon^k_0(z) = \mathbf{1}_{I_k}(z)$ is obtain using a semi-explicit Euler scheme. The advantage of the decomposition is that $\sum_k \upsilon^k(t,z) = F(z)$ for all $t>0$ and $z\in[z_{min,z_{max}}].$ Then the numerical distribution $w^i(s)$ of $Y_s$ starting from the trait $z_i$ is given by $w^i(s,y_k) = \dfrac{\upsilon^k(-s,z_i)}{F(z_i)}$ (see video and \Cref{fig lineage asex distrib}). In our simulation, we look at the particular point $z^*$ which corresponds to the dominant trait of the population, that is the mean of $F$.
		
\begin{figure}[H]
			\begin{center}
			\includegraphics[width=0.45\textwidth]{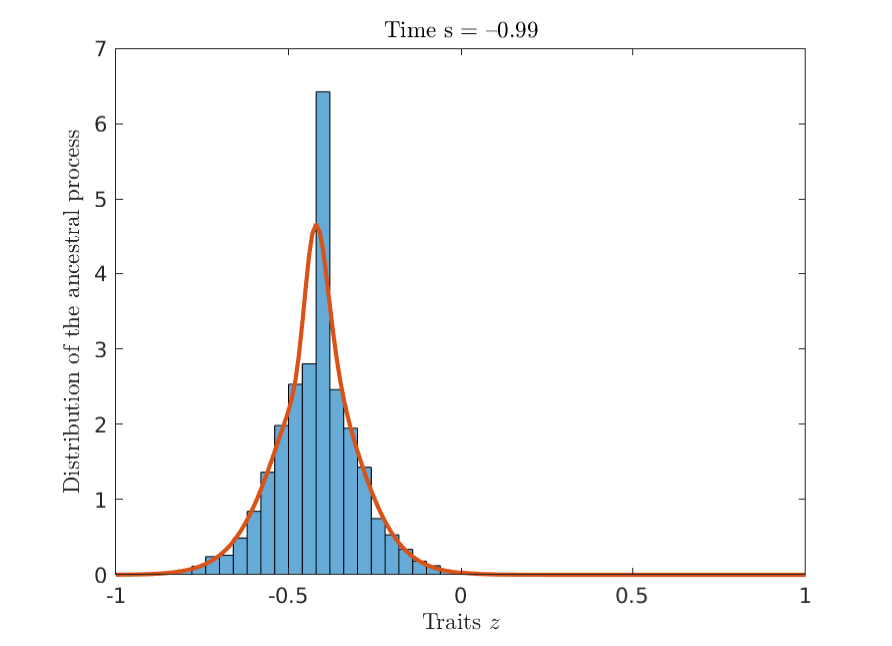}
				\caption{Ancestral lineages distribution over time of individuals with the dominant trait $z^*$. \href{https://www.lama.univ-savoie.fr/pagesmembres/garnier/Dyn_lineage_PDE_vs_IBM.mp4}{\tt  DynlineagePDEvsIBM.mp4}.}
				\label{app:movie}
			\end{center}
		\end{figure}

		As expected, we observe a good fit between the distribution of the ancestral lineage of the IBM model and the distribution of the ancestral process (see \Cref{fig lineage asex distrib} and video~\ref{app:movie}). In particular, the means coincide and each ancestral lineage of the IBM model lies in the region where the ancestral process is the most likely to be (blue region in \Cref{fig neutral frac}).


		
		\newpage
		
				\enlargethispage{-8\baselineskip}

		\subsection{Deterministic approximation of the ancestral lineage}
				\begin{figure}[h]
			\begin{center}
				\includegraphics[scale=0.75]{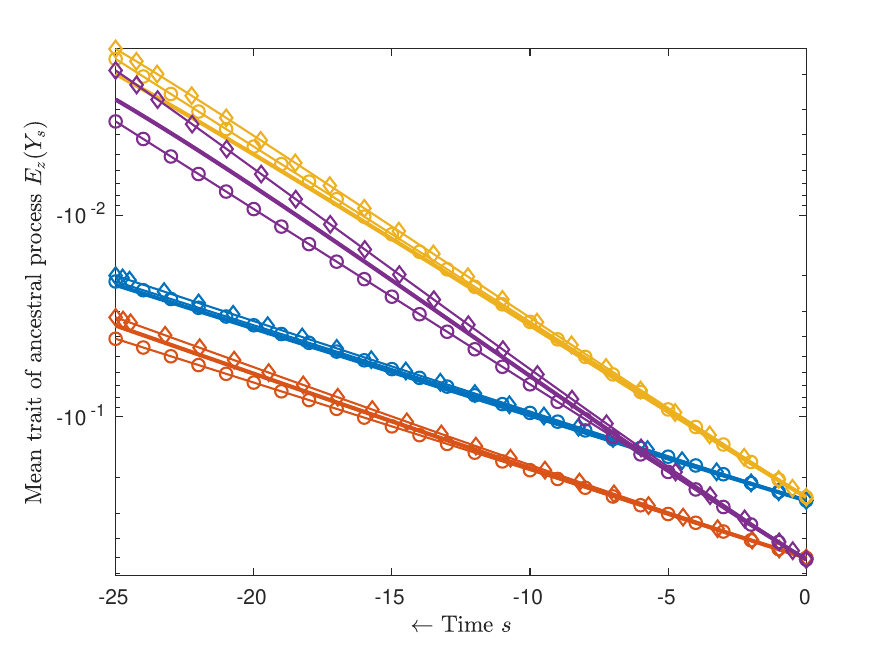}
				\caption{Mean trait of the ancestral lineage along time (plain curve), diffusive approximation (circle marked curve) and the Hamilton Jacobi approximation (diamond marked curve) in the semilog representation. The colour corresponds to different parameters: blue curve ($\sigma=0.05$, $c=0.025$), red curve ($\sigma=0.05$, $c=0.05$), orange curve ($\sigma=0.1$, $c=0.05$), purple curve ($\sigma=0.1$, $c=0.1$)  }\label{fig:mean_approx}
			\end{center}
		\end{figure}
		In this section, we aim to compare the diffusive approximation stated in \Cref{sec diffusive approx} and the Hamilton-Jacobi approximation stated in~\ref{sec HJ} with the ancestral lineage defined in~\Cref{theo edp lineages}. We compute numerically the mean of the ancestral lineage $\mathbf{E}_z(Y_s)$ using the fraction approach detailed in the above section and we compare it to the formula \cref{mean var quadex} and to the solution of the ODE~\cref{edo Gamma_bis}. Results are shown on \Cref{fig:mean_approx}.
			\end{appendices}
\end{document}